\documentclass{amsart}
\usepackage{epsfig}
\usepackage{graphics}
\usepackage{amsfonts}
\usepackage{amsmath}
\usepackage{latexsym}

\newtheorem{theorem}{Theorem}[section]
\newtheorem{lemma}{Lemma}[section]
\newtheorem{proposition}{Proposition}[section]
\newtheorem{corollary}{Corollary}[section]
\newtheorem{remark}{Remark}[section]

\numberwithin{equation}{section}

\begin{document}
\title[Multiplicative refinement paths to Brownian motion]{Bilateral Canonical
Cascades: Multiplicative Refinement Paths to Wiener's and Variant Fractional
Brownian Limits}
\author{Julien Barral}
\author{Beno\^\i t Mandelbrot}
\address{INRIA Rocquencourt, B.P. 105, 78153 Le Chesnay Cedex, France}
\email{Julien.Barral@inria.fr}
\address{Pacific Northwest National Laboratory, 222 Third Street, Cambridge,
MA \ 02142, USA}
\email{Benoit.Mandelbrot@yale.edu}
\thanks{The authors thank Pierre-Lin Pommier for kind help in the numerical simulations.}
\begin{abstract}
The original density is 1 for $t\in (0,1)$, $b$ is an integer base ($b\geq 2$%
), and $p\in (0,1)$ is a parameter. The first construction stage divides the
unit interval into $b$ subintervals and multiplies the density in each
subinterval by either $1$ or $-1$ with the respective frequencies of $\frac{1%
}{2}+\frac{p}{2}$ and $\frac{1}{2}-\frac{p}{2}$. It is shown that the
resulting density can be renormalized so that, as $n\rightarrow \infty $ ($n$
being the number of iterations) the signed measure converges in some sense
to a non-degenerate limit. If $H=1+\log _{b}$ $p>{1}/{2}$, hence $p>b^{{-1}/{%
2}}$, renormalization creates a martingale, the convergence is strong, and
the limit shares the H\"{o}lder and Hausdorff properties of the fractional
Brownian motion of exponent $H$. If $H\leq {1}/{2}$, hence $p\leq b^{{-1}/{2}%
}$, this martingale does not converge. However, a different normalization
can be applied, for $H\leq \frac{1}{2}$ to the martingale itself and for $H>%
\frac{1}{2}$ to the discrepancy between the limit and a finite
approximation. In all cases the resulting process is found
to converge weakly to the Wiener Brownian motion, independently of $H$ and
of $b$. Thus, to the usual additive paths toward Wiener measure, this
procedure adds an infinity of multiplicative paths.
\end{abstract}

\maketitle






\section{Introduction}

To motivate and clarify a new construction, this introduction compares it with others that are widely familiar. After a non-random construction has been randomized, its outcome may range
from "loosening up" slightly to changing
completely. Both possibilities, as well as intermediate ones, enter in this
paper. The point of departure is a family of non-random "cartoon" functions \cite{Maffinity} that are constructed by multiplicative interpolation. Designed as counterparts of Wiener Brownian motion \cite
{Wiener1,Wiener2} or fractional Brownian motion \cite{Kol,M2}, they have
proven to be very useful in teaching and in applications. They, in turn, are
made random in this paper, in a way that seems a "natural
inverse" but actually fails to be a straightforward step
back to the original. The fact that it reveals new interesting phenomena
suggests that the study of fractals/multifractals continues to be in large
part driven by novel special constructions with odd properties, and not only
by a general theory. Those non-random cartoons, together with a few other
examples, contradict the widely held belief that multifractal functions
(variable H\"{o}lder's $H$) are constructed by 
"multiplicative chaos" and unifractal functions (uniform H\"{o}lder's $H$), by "additive chaos".

\begin{center}
\begin{figure}
\begin{center}
{\includegraphics[width=10cm,height = 10cm]{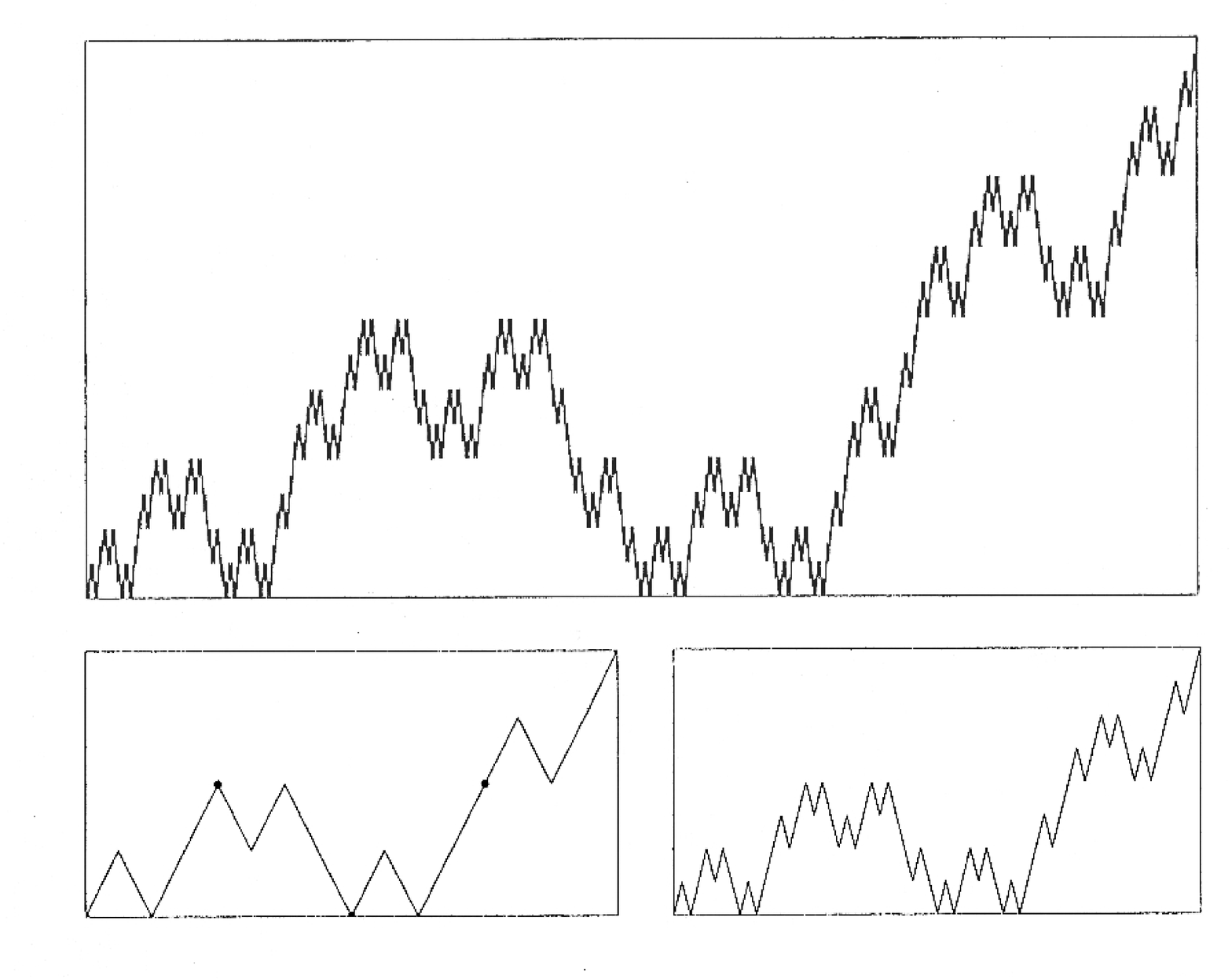}}
\vskip -.5cm
\caption{A Basic non-random cartoon of Wiener Brownian motion. Stage 2 (lower left). Stage 3 (lower right). Stage 4 (top panel).}
\label{wienerbrownianmotion}
\end{center}
\end{figure}
\end{center}

The non-random prototype described in \cite{Maffinity} is the crudest
cartoon of Wiener-Brownian motion illustrated in figure \ref
{wienerbrownianmotion}. The indicator joins the points $(0,0)$ and $(1,1)$.
The base is $b=4$ and the generator $G(t)$ is graphed by four intervals of
slope $2$ or $-2$ forming a piecewise linear continuous graph linking the
following points: $\left( 0,0\right) $, $(\frac{1}{4},\frac{1}{2})$, $(\frac{
1}{2},0)$, $(\frac{3}{4},\frac{1}{2})$, and $(1,1)$. Recursive interpolation
using this generator yields a curve characterized by the Fickian exponent $H=
\frac{1}{2}$.

A very limited randomization, described as "shuffling", moves the interval of slope $-2$ along the
abscissa from the second position to a randomly chosen position. Shuffling is
a familiar step in binomial or multinomial multifractal measures. More
interesting is the more thorough randomization introduced in \cite{M1} and called "canonical". In this context it chooses each of the four
intervals of the generator at random, independently of the others, so that
increasing and decreasing intervals have probabilities equal to their
frequencies in the original cartoon. Here $p_{+}=\frac{3}{4}$ and $p_{-}=
\frac{1}{4}$. The increment $G(1)-G(0)$ is no longer equal to $1$, but
random with the expected value $1$. As a result, the construction is no
longer a recursive interpolation and can be called a recursive refinement.

A more general construction of a non-random cartoon has an arbitrary base $
b>3$ and a continuous piecewise linear generator made of $b$ intervals of
slope $\frac{b}{c}$, where $c\le b$ is a second integer base. In that case,
defining $H$ by $c=b^{H}$ and $p$ as $p=b^{H-1}$, the frequencies $f_{+}$
and $f_{-}$ are $f_{+}=\frac{1}{2}+\frac{1}{2}\frac{c}{b}=\frac{1}{2}+\frac{
b^{H-1}}{2}=\frac{1}{2}+\frac{p}{2}$ and $f_{-}=\frac{1}{2}-\frac{p}{2}$.
The limit "cartoon" is of unbounded variation and shares the H\"{o}lder and Hausdorff properties of Wiener or Fractional Brownian Motion with
the same $H$. In the Brownian motions and their cartoons, the correlations
between past and future is known to take the form $2^{2H-1}-1$. It is
positive $\in \left( 0,1\right) $ in the persistent case $H>\frac{1}{2}$,
and negative $\in \left( -\frac{1}{2},0\right) $ in the antipersistent case $
H<\frac{1}{2}$.

A canonical randomization of any of the Brownian cartoons can now be
described. A first step consists in making all those frequencies into
probabilities. A second step consists in eliminating various constraints on $
H$ and $p$ that are due to their origin in cartoons. Both the Wiener and
Fractional Brownian motions non-random cartoons require that $
0<H<1 $ and that $H$ be a ratio of logarithms of integers, of the form $
\frac{\log c}{\log b}$. We shall allow $p$ to vary from $0$ to $1$, which
implies $-\infty <H<1$, and leave $c$ unrestricted, allowing it even to be
smaller than 1.

This paper's object is to describe the limits of the functions $\left\{
B^H_{n}\right\} _{n\geq 1}$ generated by this procedure. 
The derivatives of the functions $\left\{
B^H_{n}\right\} _{n\geq 1}$ (in the sense of distributions) form a signed
measure-valued martingale. Moreover, $B^H_n$ is absolutely continuous, and the correlation of the derivatives of $B^H_n$ and $B^H_{n+1}$ is equal to $p$ almost everywhere. For the classic positive canonical cascades \cite{M1}, those martingales converge strongly to a limit.
But that limit can degenerate to 0, and, if so, no normalization
yielding a non zero limit is known. 

\begin{center}
\begin{figure}
\begin{center}
{\includegraphics[width=10cm,height = 10cm]{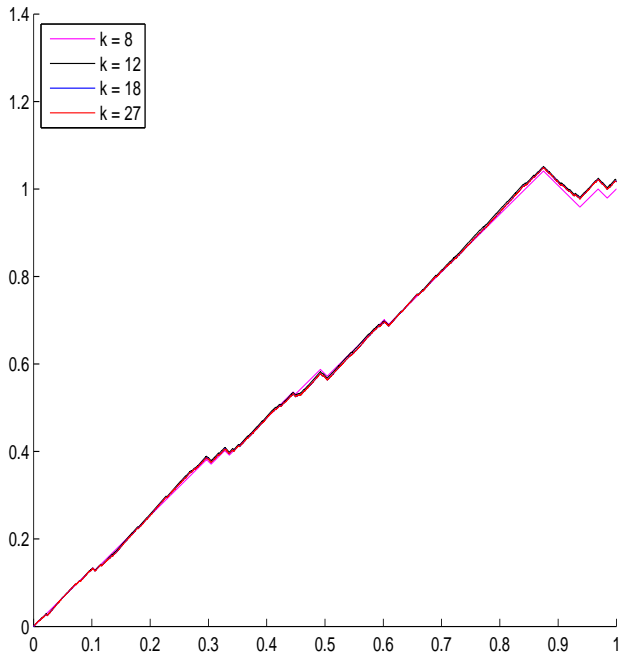}}
\vskip -1cm
\caption{$B^H_k$ for $k=8,\ 12,\ 18,\ 27$ in the case $b=2$ and $H=0.95$: Fast strong convergence.}
\label{H=.95}
\end{center}
\end{figure}
\end{center}

\begin{center}
\begin{figure}
\begin{center}

{\includegraphics[width=10cm,height = 10cm]{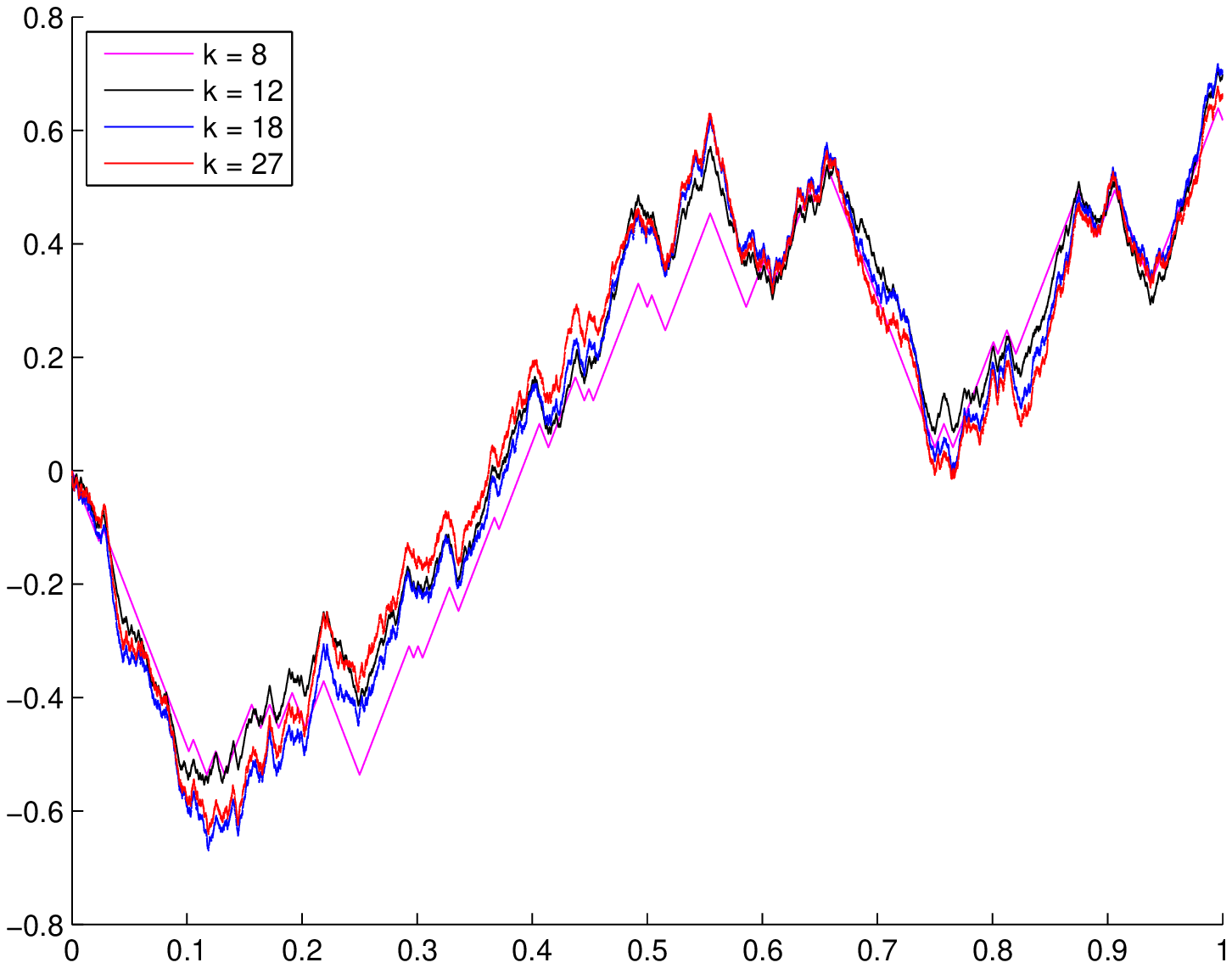}}
\vskip -.5cm
\caption{$B^H_k$ for $k=8,\ 12,\ 18,\ 27$ in the case $b=2$ and $H=0.7$. Strong convergence.}
\label{H=.7}
\end{center}
\end{figure}
\end{center}

\begin{center}
\begin{figure}
\begin{center}
{\includegraphics[width=10cm,height = 10cm]{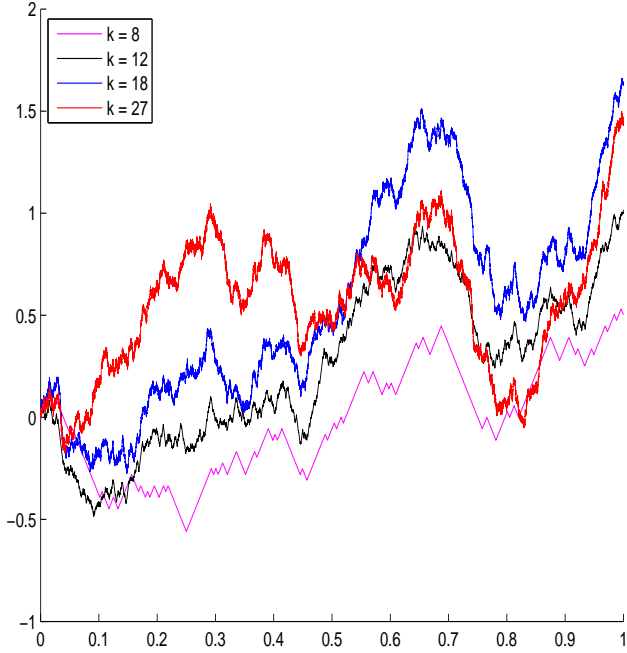}}
\vskip -1cm
\caption{$B^H_k/\sigma_{1/2} \sqrt{k}$ for $k=8,\ 12,\ 18,\ 27$ in the case  $b=2$ and $H=0.5$: Convergence in distribution to the Wierner Brownian motion.}
\label{H=.5}
\end{center}
\end{figure}
\end{center}

\begin{center}
\begin{figure}
\begin{center}
{\includegraphics[width=10cm,height = 10cm]{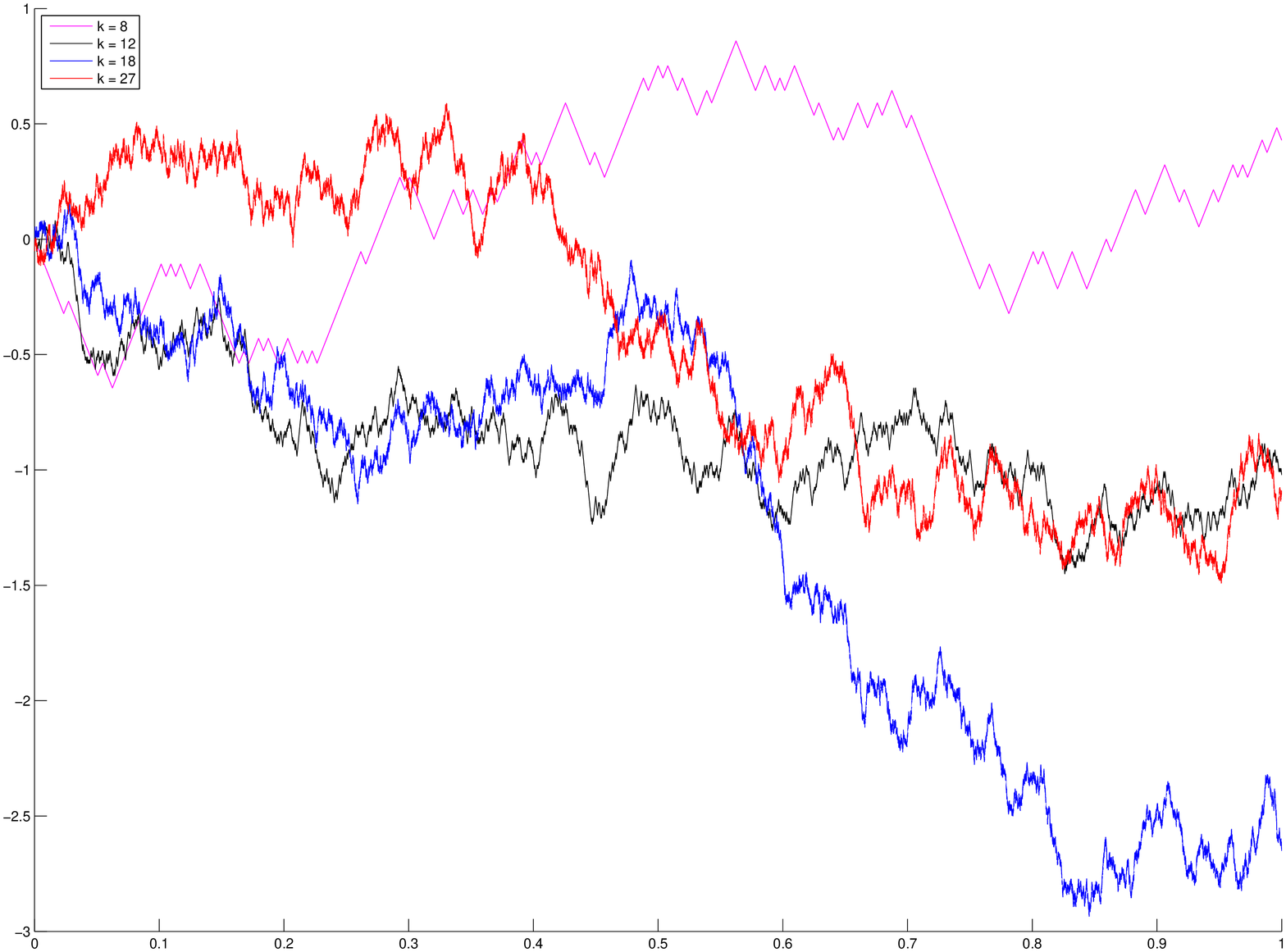}}
\vskip -1cm
\caption{$B^H_k/\sigma_H b^{k(1/2-H)}$ for $k=8,\ 12,\ 18,\ 27$ in the case  $b=2$ and $H=0.25$: Convergence in distribution to the Wierner Brownian motion.}
\label{H=.25}
\end{center}
\end{figure}
\end{center}

\begin{center}
\begin{figure}
\begin{center}
{\includegraphics[width=10cm,height = 10cm]{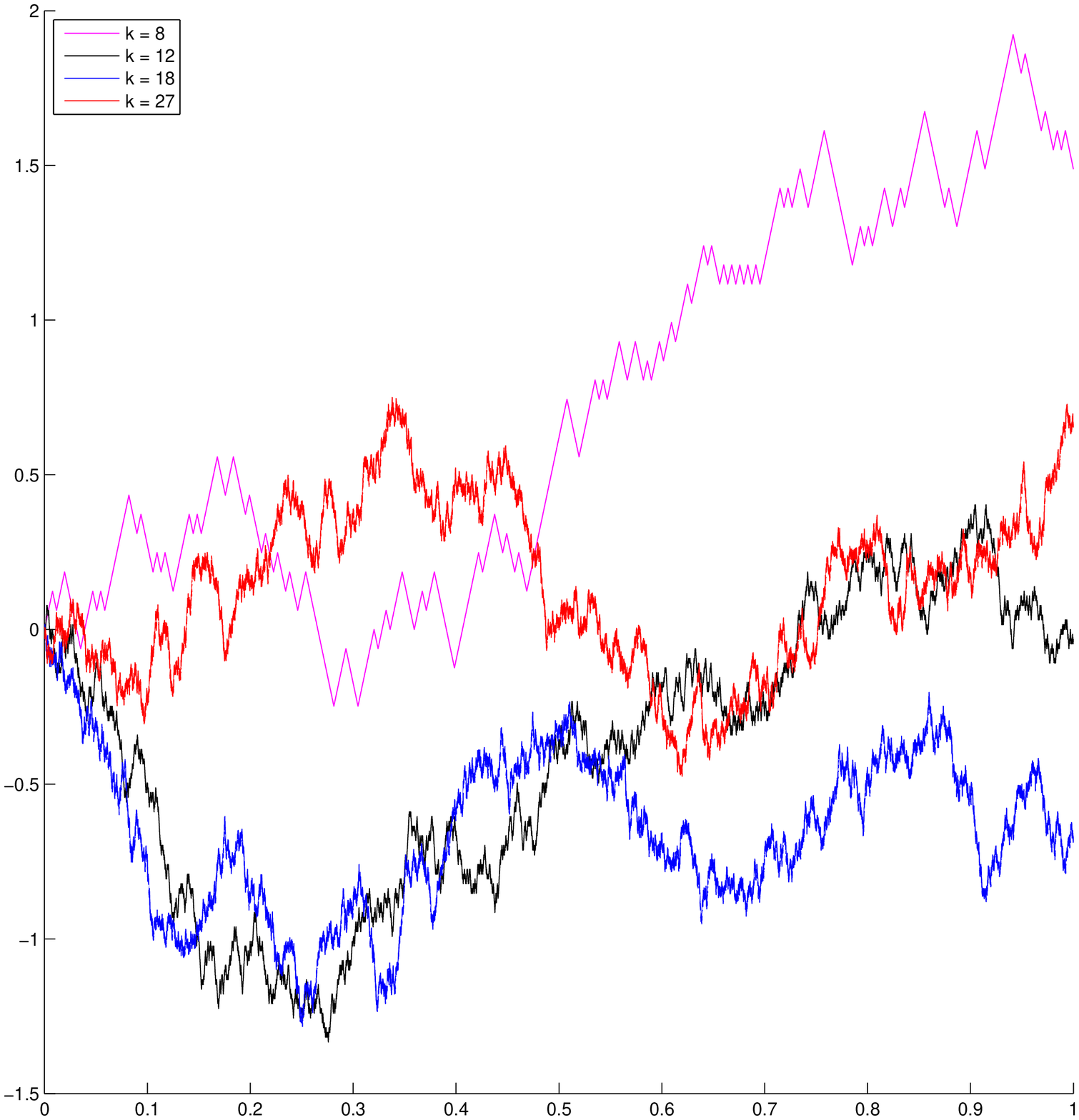}}
 \caption{$B^H_k/\sigma_H b^{k(1/2-H)}$ for $k=8,\ 12,\ 18,\ 27$ in the case $b=2$ and $H=-2$: Convergence in distribution to the Wierner Brownian motion.}
\label{H=-2}
\end{center}
\end{figure}
\end{center}

For the bilateral canonical cascades
considered in this paper, the situation will be shown to be altogether
different, following a pattern first observed in exploratory simulations. The persistent case $\frac{1}{2}<H<1$ behaves as expected, as
illustrated in Figures \ref{H=.95} and \ref{H=.7}: the martingale
converges and has a strong limit, namely, a non-Gaussian variant process sharing the H\"{o}lder and Hausdorff properties of the persistent Fractional Brownian
Motion of exponent $H$ (Theorem~\ref{th2}). The antipersistent case $-\infty <H\leq \frac{1}{2}$
, to the contrary, defies facile extrapolation. As illustrated in Figures \ref
{H=.5}, \ref{H=.25} and \ref{H=-2}, the
martingale does not converge to zero but oscillates increasingly wildly.
However, there does exist an alternative normalization that yields a
nondegenerate weak limit for $n\rightarrow \infty $, namely, Wiener Brownian
Motion. The fact that the exponent $H$ no longer affects the limit is a surprising form of
what physicists call "universality", a
phenomenon that recalls the Gaussian central limit theorem. It expresses
that the rules of the cascade are destroyed in the limit, leaving only an
accumulation of noise. Our result provides new functional central limit
theorems. Moreover, the normalization factor in the special case $H=1/2$ is
atypical (see Theorem~\ref{th1'} and Corollary~\ref{FTLC'}). Since $H$ is no
longer a H\"{o}lder, negative values of $H$ create no paradox whatsoever.

To ensure a long-range power-law correlation function, exquisite long range order must be present in Fractional Brownian
cartoons. Under canonical randomization this order is robust in the case of persistence. But it is
non robust and destroyed in the case of antipersistence, with a clear
critical point in the Wiener Brownian case. 

\medskip

This is novel but recalls an
observation concerning the Cauchy-L\'{e}vy stable exponent $\alpha $: it is
constrained to $0\leq \alpha <2$, that is, $H>\frac{1}{2}$. The non-random cartoon of one such process \cite{Mfinance} has a generator joining
the points $(0,0)$, $\left( \frac{1}{2},p\right) $, $\left( \frac{1}{2}
,1-p\right) $, $\left( 1,1\right) $. For all $p\in \left( 0,1\right) $, this
can be interpolated into a discontinuous function in which the
discontinuities $\Delta $ have a distribution of the form $\Pr \left\{
\Delta >\delta \right\} \sim \delta ^{-\alpha }$ with $\alpha =\frac{-1}{
\log _{2}p}$. This exponent can range as $\alpha \in \left( 0,\infty \right)
$. But let us, after $k$ stages, change the number of discontinuities of
given size and the number of continuous steps. Rather than fixed, let them be 
random Poissonian with the same expectation. When $\alpha <2$, the process
converges to a stable one, but when $\alpha >2$, it explodes.

In any event, for $H<\frac{1}{2}$ the integral of the covariance of the
Fractional Brownian Motion vanishes. This demands very special correlation
properties that are easily destroyed by diverse manipulations. "Instability"
also characterizes for $H<\frac{1}{2}$ the limits of expressions of the form
$\Sigma G[\Delta B_{H}(t)]$, where $\Delta B_{H}(t)$ is the increment of a
Fractional Brownian Motion over $[t,t+1]$ and $G$ is a strongly non linear
transform \cite{Taqqu}. In this context, the fact that Kolmogorov's
turbulence takes on the unstable value $H=\frac{1}{3}$ may reward a close
look.

In the preceeding first-approximation results, we perceived an analogy with
the usual sums of iid random variables $X$ of finite variance. The strong
law of large numbers tells us that the sample average, a normalized sum of $
n $ variables $X_{k}$, strongly converges to $\mathbb{E}(X)$. Then the central limit
theorem tells us that the discrepancy between $\mathbb{E}(X)$ and the $n$th normalized
sum can be subjected to a different normalization --- division
by $\sqrt{n}$ --- and after that converges
weakly to a Wiener Brownian motion. When $H\neq 1/2$, this situation generalizes to our canonical cascades,  but with some major changes. Here, $\mathbb{E}(X)$
is replaced for $H>\frac{1}{2}$ by a variant fractional Brownian motion and
for $H<\frac{1}{2}$ by $0$. The rate of convergence for the remainder depends on $H$. It takes
the same analytic form $b^{n\left( \frac{1}{2}-H\right) }$ for all $H\neq 1/2$ but plays different roles: For $H<\frac{1
}{2}$ it compensates for boundless growth, and for $H>1/2$, for a decrease to $0$ (Theorems~\ref{th1} and \ref{th3}).

Our construction possesses a natural extension to the case $H=-\infty$ if we consider the processes $B^H_n/b^{-nH}$ and let $H$ tend to $-\infty$ for every $n\ge 1$. In this case, $p=0$ and $k$
recursions yield $b^{k}$ values of a random walk with symmetric correlated increments taking values in $\{-1,1\}$. After
division by $b^{\frac{k}{2}}$, there is a limit in distribution for the associated piecewise linear function, namely, the Wiener Brownian Motion (see Corollary~\ref{FTLC}). This is an unexpected extension of the usual result known for the classical random
walk obtained by coin tossing. Here, convergence is as "weak" as can be, since the
terms in the sequence are statistically independent.

Understanding bilateral cascades is helped by a step that has been
fruitful since the earliest canonical cascades \cite{M1}. It consists in
keeping $p$ constant, replacing the interval $[0,1]$ by the cube $[0,1]^{E}$, and varying the Euclidean dimension $E$ from a large value down. In all
cascade constructions, the proper distance is not Euclidean but ultrametric.
Hence, the nondegenerate versus degenerate alternative requires no new
argument: it proceeds just as on a linear grid of base $b^{E}$. The critical
$H=\frac{1}{2}$ now corresponds to $p_{crit}=b^{\frac{-E}{2}}$, which $
\searrow 0$ as $E\nearrow \infty $. In a high-dimensional space, a cascade
with the given $p$ yields a variant fractional Wiener signed measure but the
intersections of that measure by subspaces of small\textbf{\ }$E$
degenerates to an infinitesimal Wiener measure (this extension to higher
dimensions will be studied in a further work). Classically, this is also the
case in birth and death cascades with multiplier values 1 and 0. The novelty
present in the bilateral case in that the term "degenerate" takes a
different meaning.

The martingales considered in this paper are the very simplest special case of the following more general construction. As for positive canonical cascades, given an integer $b\ge 2$, the recursive process consists in associating with each $b$-adic subinterval $J$ of $[0,1]$ a random weight $W_J$ so that these weights are i.i.d with a random variable $W$ and $\mathbb{E}(W)$ is defined and equal to $1/b$.  Then, one gets a sequence of random piecewise functions $(F_n)_{n\ge 1}$ by imposing that $F_n(0)=0$ and that the increment of $F_n$ over the interval $J$ of the $n^{\mbox{{\small th}}}$ generation is equal to the product $W_{J_1}W_{J_2}\cdots W_{J_n}$, where $J_k$ is the $b$-adic interval of generation $k$ containing $J(=J_n)$. Observe that this construction falls in the category of infinite products of functions \cite{Fan,BM}.  The family $\{F_n\}_{n\ge 1}$ forms a continuous functions-valued martingale. A sufficient condition for the sequence $F_n$ to converge almost surely uniformly is that  the function $\tau_W(p)=q-1-\log_b\mathbb{E}(|W|^p)$ takes a positive value at some $p\in (1,2]$. In the simplest case studied in this paper, $W$ belongs to $\{-b^{-H},b^{-H}\}$, and the critical value $H=1/2$ separates the domain $H\le 1/2$ for which $\tau_W(p)\le 0$ over $[1,2]$ and the domain $H\in (1/2,1]$ for which we always have $\tau_W(2)>0$. In the general case, when $\tau_W((1,2])\not\subset (-\infty,0]$, the limit of the signed canonical cascade is not a unifractal but a multifractal function -- to be studied in a further work.

The rest of the paper is organized as follows. This section ends with definitions and notations used in the sequel. Also, the processes studied in this paper are more formally defined than in the previous paragraphs. Section~\ref{statements} and~\ref{statements2} provide our main results for the cases  $p\le b^{-1/2}$ (i.e. $H\le 1/2$) and $p>b^{-1/2}$ (i.e. $H>1/2$) respectively. The next three sections are devoted to the proofs of our main results.  


\section{Construction of the martingale}
\subsection{Definitions and notations}\label{def}
$\ $

Let $b\geq 2$ be an integer.

For $n\ge 0$ let $\Sigma_n=\{0,\dots,b-1\}^{n}$, where $\Sigma_0$ contains
only the empty word denoted by $\emptyset$. Also let $\Sigma^*=\bigcup_{n\ge
0}\Sigma_n$ and $\Sigma=\{0,\dots,b-1\}^{\mathbb{N}^*}$. The concatenation
operation from $\Sigma^*\times (\Sigma^*\bigcup \Sigma)$ to $
(\Sigma^*\bigcup \Sigma)$ is denoted $\cdot$.

For $x\in \Sigma $ and $n\geq 1$, let $x|n$ be the projection of $x$ on $
\Sigma _{n}$ and $x|\infty=x$. Then for $n\geq 1$ and $w\in \Sigma _{n}$, we set $[w]=\{x\in
\Sigma :x|n=w\}$. Given two words of infinite length $x,y\in \Sigma $, one
defines $x\wedge y$ as $x|n_{0}$, where $n_{0}=\sup \{n\geq
1:x|n=y|n\} $. Adopt the convention that $\inf \emptyset =0$ and $x|0$
is the empty word $\emptyset $.

The length of any element $w$ of $\Sigma_n$ is equal to $n$ and is denoted
by $|w|$.

Denote by $\pi $ the mapping $x\in \Sigma \mapsto \sum_{k=1}^{\infty
}x_{k}b^{-k}$.

If $w\in \Sigma^*$, $t_w$ stands for the number $\sum_{k=1}^{|w|}w_kb^{-k}$
and $I_w$ stands for the closed $b$-adic interval $\pi ([w])$.

\smallskip

For $n\geq 0$ denote by $\mathcal{D}_{n}$ the set of $b$-adic numbers of the 
$n^{\mbox{{\small th}}}$ generation in $[0,1]$. Also denote by $\mathcal{D}$ the set
of all $b$-adic numbers of $[0,1]$.

\smallskip

Denote by $\mathcal{C}([0,1])$ the space of real valued continuous
functions on $[0,1]$. Then, for $\alpha\in (0,1]$, $\mathcal{C}
^\alpha([0,1]) $ stands for the subspace of $\mathcal{C}([0,1])$ whose
elements are uniformly $\alpha$-H\"older continuous, i.e. $f\in \mathcal{C}^{\alpha}
([0,1])$ if and only if there exists $C>0$ such that $|f(t)-f(s)|\le
C|t-s|^\alpha$ for all $t,s\in [0,1]$.

If $f\in \mathcal{C}([0,1])$, denote its modulus of continuity by $\omega
(f,\cdot )$ (for $\delta \in \lbrack 0,1]$, $\omega (f,\delta )=\sup_{t,s\in
\lbrack 0,1],|t-s|\leq \delta }|f(t)-f(s)|$).

Recall that the pointwise H\"older exponent of $f$ at $t_0\in [0,1]$ is
defined by 
\begin{equation*}
h_f(t_0)=\sup \left \{\alpha\ge 0: \ \exists \ P\in\mathbb{R}[X],\
\sup_{t\in [0,1]\setminus \{t_0\}}\frac{ |f(t)-f(t_0)-P(t)|}{|t-t_0|^\alpha}
<\infty\right \}.
\end{equation*}

If $I$ is a subinterval of $[0,1]$, $\Delta f(I)$ stands for $
|f(\sup(I))-f(\inf(I))|$.

\subsection{A construction of a recursive canonical cascade with values $\pm
1$.}
$\ $

\noindent 
Let $(\Omega ,\mathcal{A},\mathbb{P})$ be the probability space on which the
random variables in the sequel are defined. If $Y$ is a random variable, we
shall denote by $\mathcal{L}(Y)$ its probability distribution.

For $0\le k\le b-1$ let $S_k(t)=b^{-1}t+k$.

If $H\in \lbrack -\infty ,1]$, define the probability measure $\pi _{b,H}=p_{b,H}^{+}\,\delta
_{1}+p_{b,H}^{-}\,\delta _{-1}$, where

\begin{equation*}
p_{b,H}^+=\frac{1+b^{H-1}}{2}\quad \mbox{and} \quad p_{b,H}^-=1-p_{b,H}^+,
\end{equation*}
with the convention $b^{-\infty}=0$.

For all $H\in \lbrack -\infty ,1]$, let $(\epsilon ^{H}(w))_{w\in \Sigma
^{\ast }}$ be a sequence of mutually independent random variables of common
probability distribution $\pi _{b,H}$. When $H$ is fixed in the sequel,
sometimes we simply write $\epsilon (w)$ for $\epsilon ^{H}(w)$.

\medskip

If $u\in \lbrack 0,1]\setminus \mathcal{D}$, we identify $u$ with the unique
element $\widetilde{u}\in \Sigma $ such that $u=\pi (\widetilde{u})$. Then,
for every $n\geq 1$, if $H\in (-\infty ,1]$ we consider on $[0,1]$ the continuous
piecewise linear map $B^H_n$ over the $b$-adic intervals of the $n^{\mbox{\small th}}$ generation 
such that $B^H_n(0)=0$ and for $w \in \Sigma_n$ the increment of $B^H_n$ over $I_w$ is equal to $\epsilon^H(w|1) \cdots \epsilon^H (w|n) b^{-nH}$, i.e.
\begin{equation*}
B^H_n (t)=b^{-nH}\int_0^t b^n \epsilon^H(u|1) \cdots \epsilon^H (u|n)\, du,
\end{equation*}
We leave the reader verify that the sequence $(B^H_n)_{n\ge 1}$ is a $\mathcal{C}([0,1])$-valued martingale
with respect to the filtration $\big (\sigma(\epsilon (w): \ w\in \Sigma_n)
\big )_{n\ge 1}$.

More generally for all $w\in \Sigma^*$ let 
\begin{equation*}
B^H_n(w) (t)=b^{-nH}\int_0^t b^n \epsilon^H(w\cdot u|1) \cdots \epsilon^H
(w\cdot u|n)\, du.
\end{equation*}
Of course, $B^H_n(\emptyset)=B^H_n$ almost surely.

For all $0\le k\le b-1$ and $t\in I_k=[k/b,(k+1)/b]$ we have the relation 
\begin{equation}  \label{fonc}
B^H_n(t)-B^H_n(k/b)= \epsilon^H(k) b^{-H} B_{n-1}^H(k)\left (S_{k}^{-1}
(t)\right )
\end{equation}
and more generally for all $w\in\Sigma^*$, $t\in I_w$ and $n\ge |w|$

\begin{equation}  \label{fonc0}
B^H_n(t)-B^H_n(t_w)= \epsilon^H (w|1) \cdots \epsilon^H (w|n) b^{-nH}
B_{n-|w|}^H(w)\left (S_{w_n}^{-1}\circ \cdots \circ S_{w_1}^{-1} (t)\right ).
\end{equation}

For $n\geq 0$ and $w\in \Sigma ^{\ast }$ we denote by $Z_{n}(w)$ the random
variable $B_{n}^{H}(w)(1)$, with the convention $B_{0}^{H}(w)(1)=1$. When $
w=\emptyset $, we simply write $Z_{n}$ for $Z_{n}(w)$. The relation (\ref
{fonc}) yields for every $n\geq 1$ 
\begin{equation} \label{fonc2}
Z_{n}=\sum_{k=0}^{b-1}b^{-H}\epsilon _{k}Z_{n-1}(k), 
\end{equation}
where the random variables $\epsilon _{0},\dots ,\epsilon
_{b-1},Z_{n-1}(0),\dots ,Z_{n-1}(b-1)$ are mutually independent. Moreover, $
\mathcal{L}(\epsilon _{k})=\pi _{b,H}$ and $\mathcal{L}(Z_{n-1}(k))=\mathcal{
L}(Z_{n-1})$ for all $0\leq k\leq b-1$. Relation (\ref{fonc2}), which will
be useful in the sequel, is familiar from the positive cascade case \cite{M1}.

Finally, for $H\in [-\infty,1/2]$ let 
\begin{equation*}
\sigma _{H}=
\begin{cases}
1 & \mbox {if $H=-\infty$} \\ 
\displaystyle\sqrt{\frac{b-1}{b^{2-2H}-b}+1} & 
\mbox {if $H\in
(-\infty,1/2)$} \\ 
\displaystyle\sqrt{\frac{b-1}{b}} & \mbox {if $H=1/2$}
\end{cases}.
\end{equation*}

\section{Weak convergence of the normalized martingale $B^H_n$ to Wiener Brownian motion, independently of $b$ and $H$ in the anti-persistant case H $\leq $ 1/2}~\label{statements}

\begin{theorem}
\label{th1} Let $H\in (-\infty,1/2)$. The sequence $\left (\mathcal{L}\big(B_{n}^{H}/\sigma_H  b^{n(1/2-H)}\big)\right )_{n\geq 1}
$ converges weakly to the Wiener measure as  $n$ goes to $\infty $.
\end{theorem}

\begin{theorem}\label{th1'}
The sequence $\left (\mathcal{L}\big(B_{n}^{1/2}/\sigma_{1/2}\sqrt{n}\big)\right )_{n\geq 1}
$ converges weakly to the Wiener measure as  $n$ goes to $\infty $.
\end{theorem}

\begin{remark}\label{dege}
(1) When $H<1/2$, Theorem~\ref{th1} implies that with probability 1, $
\limsup_{n\rightarrow \infty }\frac{\Vert B_{n}^{H}\Vert _{\infty }}{b^{n(1/2-H)}}
>0$. Thus the martingale $B^H_n$ neither strongly converges to a non trivial limit in $\mathcal{C}([0,1])$ nor to 0. This fact deserves to be called degeneracy. This shows a strong difference with positive canonical cascades for which degeneracy means uniform convergence to 0 (\cite{M1,KP}). The same remarks hold when $H=1/2$.

(2) For $H\le 1/2$ define $X^H_n=B^H_n(1)/\sigma_H  b^{n(1/2-H)}$ if $H< 1/2$ and $X^H_n=B_{n}^{1/2}(1)/\sigma_{1/2}\sqrt{n}$ otherwise. The reader can check that $X^H_n(1)$ is not a Cauchy sequence in $L^2$ while the $L^2$ norm of $X^H_n(1)$ converges to 1 as $n$ goes to $\infty$. This implies that  $X_n^H(1)$ cannot converge almost surely to a standard normal random variable. Consequently, Theorem~\ref{th1} and \ref{th1'} cannot be strengthen into results of almost sure convergence.

\end{remark}

\section{Restatement of Theorems~\ref{th1} and ~\ref{th1'} as functional CLT with atypical normalization when $H=1/2$}\label{TLC}$\ $

If $H\in [-\infty,1/2]$, $n\geq 1$ and $0\leq k<b^{n}$ and $w$ is the unique element of $\Sigma
_{n}$ such that $t_{w}=kb^{-n}$ let 
\begin{equation*}
\xi _{k}^{(n,H)}=\prod_{j=1}^{n}\epsilon ^{H}(w|j).
\end{equation*}
For a given $n\ge 1$, the random variables $\xi
_{k}^{(n,H)}$, $0\le k<b^n$, are identically distributed, and they take values in $\{-1,1\}$.

Also, consider the random walk $\big (S_{r}^{(n,H)}\big)_{0\le r<b^n}$ defined by
\begin{equation*}
S_{r}^{(n,H)}=\sum_{k=0}^{r-1}\xi _{k}^{(n,H)}
\end{equation*}
(with the convention $S_{-1}^{(n,H)}=0$). 

\begin{corollary}[Functional central limit theorem]
\label{FTLC} Let $H\in [-\infty,1/2)$ and for $n\ge 1$ and $t\in [0,1]$ define  $X_{n}^{H}(t)=\displaystyle\frac{1}{\sigma _{H}\sqrt{b^n}}\left [
S_{[b^{n}t]}^{(n,H)}+(b^{n}t-[b^{n}t])\xi _{\lbrack b^{n}t]}^{(n,H)} \right ]$. The sequence $\mathcal{L}(
X^H_n)_{n\ge 1}$ converges weakly to the Wiener measure as $n$ tends to $
\infty$.
\end{corollary}

\begin{corollary}[Functional central limit theorem]
\label{FTLC'} For $n\ge 1$ and $t\in [0,1]$ define  $X_{n}^{1/2}(t)= \displaystyle\frac{1}{\sigma _{1/2}\sqrt{nb^{n}}}
\left [S_{[b^{n}t]}^{(n,1/2)}+(b^{n}t-[b^{n}t])\xi _{\lbrack b^{n}t]}^{(n,1/2)} \right ]$. The sequence $\mathcal{L}(
X^{1/2}_n)_{n\ge 1}$ converges weakly to the Wiener measure as $n$ tends to $
\infty$.
\end{corollary}

\begin{remark}
(1) When $H<1/2$, the stochastic process $X_{n}^{H}$ takes formally the
same form as a the processes considered in central limit theorems for weakly
dependent sequences (see \cite{Bil}, Ch. 19 or \cite{Doukhan}). The main difference is that in
the process we consider the random variables $\xi _{k}^{(n,H)}$ are highly
correlated. Nevertheless, the same asymptotic behavior (weak convergence to
the Wiener measure) holds.

(2) In the case $H=1/2$, the normalizing factor takes an untypical form.

\end{remark}

\begin{remark}
When $H=-\infty$, it is not
difficult to verify that the conclusions of Propositions~\ref{prop1}--\ref
{prop3} hold for $(X^{-\infty}_n)_{n\ge 1}$ instead of the sequence $(B^H_n/a^H_n)_{n\ge
1}$ considered in Section~\ref{proofth1} by using the same approach. Consequently, the proof
is left to the reader in this case.
\end{remark}

\begin{remark}
We mention that functional central limit theorems associated with positive canonical cascades have been established in \cite{LRR} in a very different spirit. There, a square integrable random weight $W$ is fixed which generates a canonical multiplicative cascades in base $b$ and its associated sequence of increasing functions $F^{(b)}_n$ on $[0,1]$ converging to a function $F^{(b)}_\infty$. For each $ n\in \mathbb{N}^*\cup\{\infty\}$ the authors establish a functional central limit theorem for $\big (F_n^{(b)}(t)-t\big )/\sqrt {b}$ as the basis $b$ tends to $\infty$: $\big (F_n^{(b)}(t)-t\big )/\sqrt {b}$ converges in law to a multiple of the Brownian motion. As letting $b$ tend to $\infty$ weakens the correlations between the increments of $F^{(b)}$, the existence of such a weak limit is natural. 
\end{remark}

\section{Strong convergence of the martingale in the persistent case ${1}/{2}<H\le 1$}~\label{statements2}

\begin{theorem}\label{th2}
Suppose that $H\in (1/2,1]$. The sequence $(B_{n}^{H})_{n\geq 1}$
is a martingale that converges almost surely and in $L^{2}$ norm to a
continuous function $B^{H}$. Moreover, with probability 1,

\begin{enumerate}
\item $B^{H}$ belongs to $\bigcap_{H'<H}\mathcal{C}^{H'}([0,1])$ and it has everywhere a pointwise H\"{o}lder exponent equal to $H$.

\item The Hausdorff and box dimensions of the graph of $B^{H}$ is $2-H$.
\end{enumerate}
\end{theorem}

\begin{remark}
The limit process $B^H$ is not Gaussian since a computation shows that the third moment of the centered random variable $B^H(1)-1$ does not vanish. 

Notice that the case $H=1$ yields the deterministic function $B^1(t)=t$.
\end{remark}

\section{Functional CLT associated with the strong convergence case $1/2<H<1$}$\ $

It will be shown in Section~\ref{strong} that $\mathbb{E}(B^H(1)^2)<\infty$ if $H>1/2$. Consequently the number $\sigma_H=\sqrt{\mathbb{E}(B^H(1)^2)-1}$ is positive and finite when $H\in (1/2,1)$.

The following Theorems~\ref{th3} and \ref{th4} must be viewed as  couterparts of Theorems~\ref{th1} and Corollary~\ref{FTLC}.

\begin{theorem}\label{th3}
Let $H\in (1/2,1)$.  The sequence $\left (\big(B^H-B^H_n\big )/\sigma_H b^{n(1/2-H)}\right )_{n\ge 1}$ converges weakly to the Wiener measure as $n$ tends to $\infty$.
\end{theorem}

If $H\in (1/2,1)$, for every $w\in \Sigma^*$ denote by $B^H(w)$ the almost sure limit of $B^H_n(w)$. Also,  if $n\geq 1$ and $0\leq k<b^{n}$ and $w$ is the unique element of $\Sigma
_{n}$ such that $t_{w}=kb^{-n}$ let 
$$
\widetilde \xi _{k}^{(n,H)}=\big (B^H(w)(1)-1\big )\prod_{j=1}^{n}\epsilon ^{H}(w|j).
$$
Then define $\displaystyle S_{p}^{(n,H)}=\sum_{k=0}^{p-1}\widetilde \xi _{k}^{(n,H)}$ for $0\le p<b^n$ and finally consider on $[0,1]$ the piecewise linear function
$$
X^H_n(t)=\frac{1}{\sigma _{H}\sqrt{b^{n}}}\left [
S_{[b^{n}t]}^{(n,H)}+(b^{n}t-[b^{n}t])\widetilde \xi _{\lbrack b^{n}t]}^{(n,H)}\right ].
$$
\begin{theorem}\label{th4}
Let $H\in (1/2,1)$. The sequence $\mathcal{L}(
X^H_n)_{n\ge 1}$ converges weakly to the Wiener measure as $n$ tends to $
\infty$.
\end{theorem}

\begin{remark}
Theorem~\ref{th3}  implies that $(B^H-B^H_n)(1)/\sigma_H b^{n(1/2-H)}$ converges in law to a $\mathcal{N}(0,1)$ law.  This result is of the same nature as Proposition 4.1 in \cite{OW} which deals with central limits theorems associated with non negative canonical cascades. The technique used in \cite{OW} would work to establish the convergence of   $\mathcal{L}\big ((B^H-B^H_n )(1)/\sigma_H b^{n(1/2-H)}\big )$. It uses Lindeberg's theorem, while we exploit the functional equation (\ref{fonc2}).
\end{remark}

\section{Proof of Theorem~\ref{th1} and \ref{th1'} and their corollaries concerning the case $H\leq \frac{1}{2}$, i.e., $
p\leq b^{-1/2}$}\label{proofth1}
Theorems~\ref{th1} and \ref{th1'} follow from the next three propositions. In fact we are going to show that the sequence $\left (\mathcal{L}\big (B^H_n/a^H_n\big )\right )_{n\ge
1}$ converges weakly to the Wiener measure, where for $n\ge 1$ $\displaystyle a_{n}^{H}=
\displaystyle\sqrt{\left( \frac{b-1}{b}+b^{1-2H}-1\right) \frac{b^{n(1-2H)}-1
}{b^{1-2H}-1}}$ if $H<1/2$ and $\displaystyle a^{H}_n=\sqrt{\frac{b-1}{b}n}$ if $H=1/2$ (observe that by the L'Hospital rule, $
a_{n}^{H}$ converges to $a_{n}^{1/2}$ as $H\nearrow
1/2 $). It is easily seen that this will imply Theorems ~\ref{th1} and \ref{th1'} and so their corollaries Corollaries~\ref{FTLC} and \ref{FTLC'}. 

When $H<1/2$, the normalization by $a_n^H$ is more practical to use than $\sigma_H b^{n(1/2-H)}$ because it naturally appears in the $B^H_n$ asymptotic behavior's study.


\begin{proposition}
\label{prop1} Let $H\in (-\infty,1/2]$. The sequence $\big (\mathcal{L}
(B^H_n(1)/a^H_n)\big )_{n\ge 1}$ converges to $\mathcal{N}(0,1)$ as $n$ goes
to $\infty$.
\end{proposition}

\begin{proof}
Let $Y_n=B^H_n(1)/a^H_n=Z_n/a^H_n$. It is enough to show that

\begin{enumerate}
\item for every $p\ge 0$ one has the property $(\mathcal{P}_{2p})$: $
M_{2p}=\lim_{n\to\infty} \mathbb{E}(Y_n^{2p})$ exists. Moreover $M_2=1$;

\item for every $p\ge 0$ one has the property $(\mathcal{P}_{2p+1})$: $
\lim_{n\to\infty} \mathbb{E}(Y_n^{2p+1})=0$;

\item the moments of even orders obey the following induction relation valid
for $p\ge 2$:

\begin{equation*}
M_{2p}=\big (b^p-b\big )^{-1}\sum_{\substack{ 0\le
\alpha_0,\dots,\alpha_{b-1}<p  \\ \sum_{k=0}^{b-1}\alpha_k=p}} \frac{(2p)!}{
(2\alpha_0)!\cdots (2\alpha_{b-1})!} \prod_{k=0}^{b-1} M_{2\alpha_k}.
\end{equation*}
\end{enumerate}

Indeed, (1) will ensure that the sequence of probability distributions $
\mathcal{L}(Y_n)$ is tight. Moreover, it is easy to verify that a $\mathcal{N
}(0,1)$ random variable $N$ is so that its moments of even orders satisfy
the same relation as the numbers $M_{2p}$, $p\ge 1$, defined by $M_2=1$ and
the induction relation (3) (to see this, write $N$ as the sum of $b$
independent $\mathcal{N}(0,b^{-1/2})$ random variables). Consequently, since
the law $\mathcal{N}(0,1)$ is characterized by its moments, $Y_n$ must
converge in law to $\mathcal{N}(0,1)$.

\medskip

Let us establish (1), (2) and (3).

Let us take the expectation of the square of (\ref{fonc2}), by using the
fact that $\mathbb{E}(Z_n)=1$. This will explain the introduction of the
normalization factor $a_n^H$.

We have (notice that $\mathbb{E}(\epsilon_0)=b^{H-1}$)

\begin{eqnarray}  \label{m2}
\mathbb{E}(Z_n^2)&=& b^{1-2H}\mathbb{E}(Z_{n-1}^2)+ b(b-1) b^{-2H}\left (
\mathbb{E}(\epsilon_0)\right )^2 \\
&=&b^{1-2H}\mathbb{E}(Z_{n-1}^2)+\frac{b-1}{b}.
\end{eqnarray}

This yields $\mathbb{E}(Z_n^2)=(a_n^H)^2+1$. In particular, the limit $M_2$
is well defined and equals 1. Moreover, $\lim_{n\to\infty}\mathbb{E}(Y_n)=0$
since $\mathbb{E}(Z_n)=1$ and $\lim_{n\to\infty} a_n^H=\infty$.

Now let $q$ be an integer $\ge 3$. Taking the expectation of (\ref{fonc}) to
the power $q$ yields

\begin{equation}  \label{Mq}
\mathbb{E}(Z_{n+1}^q)=b^{1-qH}\mathbb{E}(\epsilon_0^q)\mathbb{E}
(Z_{n}^q)+b^{-Hq}\sum_{\substack{ 0\le \beta_0,\dots,\beta_{b-1}<q  \\ 
\sum_{k=0}^{b-1}\beta_k=q}} \gamma_{(\beta_0,\dots, \beta_{b-1})}
\prod_{k=0}^{b-1} \mathbb{E}(\epsilon_0^{\beta_k})\mathbb{E}
(Z_{n}^{\beta_k}),
\end{equation}
where $\displaystyle \gamma_{\beta_0,\dots, \beta_{b-1}} =\frac{q!}{(\beta_0)!\cdots
(\beta_{b-1})!}$.

Let us denote $\mathbb{E}(Y_n^q)$ by $M_q^{(n)}$, the set $\{(\beta_0,\dots,
\beta_{b-1})\in\mathbb{N}^b: 0\le \beta_0,\dots,\beta_{b-1}<q,\
\sum_{k=0}^{b-1}\beta_k=q\}$ by $S_q$, the ratio $\displaystyle\sqrt{\frac{n
}{n +1}}$ by $r^{(1/2)}_n$ and the ratio $\displaystyle \sqrt{\frac{b^{n(1-2H)}-1}{
b^{(n+1)(1-2H)}-1}}$ by $r_n^{(H)}$ when $H<1/2$.

Now, using that $\mathbb{E}(\epsilon_0^q)=b^{H-1}$ or $1$ depending on $q$ is
an odd or an even number, (\ref{Mq}) yields for $H\le 1/2$

\begin{equation}  \label{H}
M_q^{(n+1)}= 
\begin{cases}
\displaystyle (r^{(H)}_n)^{q}\left
(b^{-(q-1)H}M_q^{(n)}+b^{-qH}\sum_{\beta\in S_q}\gamma_{ \beta}
\prod_{k=0}^{b-1}\mathbb{E}(\epsilon_0^{\beta_k}) M_{\beta_k}^{(n)}\right )
& \mbox{if $q$ is odd}, \\ 
\displaystyle (r^{(H)}_n)^{q}\left
(b^{1-qH}M_q^{(n)}+b^{-qH}\sum_{\beta\in S_q}\gamma_{ \beta}
\prod_{k=0}^{b-1} \mathbb{E}(\epsilon_0^{\beta_k})M_{\beta_k}^{(n)}\right )
& \mbox{if $q$ is even}
\end{cases}
.
\end{equation}

We show by induction that $\big ((\mathcal{P}_{2p-1}),(\mathcal{P}_{2p})%
\big
)$ holds for $p\ge 1$, and we deduce the relation (3).

We have shown that $\big ((\mathcal{P}_{1}),(\mathcal{P}_{2})\big )$ holds.
Suppose that $\big ((\mathcal{P}_{2k-1}),(\mathcal{P}_{2k})\big )$ holds for 
$1\le k\le p-1$, with $p\ge 2$. In particular, $M_{\beta_k}^{(n)}$ goes to 0
as $n$ goes to $\infty$ if $\beta_k$ is an odd integer belonging to $[1,
2p-3]$.

Suppose $H=1/2$ and simply denote $r_n^{(1/2)}$ by $r_n$. Every element of the set $S_{2p-1}$ must contain an odd
component. Due to our induction assumption, this implies that in the
relation (\ref{H}), the term $r_n^{2p-1}b^{-(2p-1)/2}\sum_{\beta\in
S_{2p-1}}\gamma_{\beta} \prod_{k=0}^{b-1}\mathbb{E}(\epsilon_0^{\beta_k})
M_{\beta_k}^{(n)}$ in the right hand side of $M_{2p-1}^{(n+1)}$ goes to 0 at 
$\infty$. This yields 
\begin{equation*}
M^{(n+1)}_{2p-1}=r_n^{2p-1}b^{-(p-1)}M^{(n)}_{2p-1}+o(1)
\end{equation*}
as $n\to\infty$. Since $r_n^{2p-1}b^{-(p-1)}\le b^{1-p}<1$, this yields $\lim_{n\to\infty}M^{(n)}_{2p-1}=0$, that is to say $(\mathcal{P}_{2p-1})$.

Now, the same argument as above shows that in the right hand side of $
M^{(n+1)}_{2p}$, we have

\begin{equation*}
\lim_{n\to\infty} \sum_{\beta\in S_{2p}}\gamma_{ \beta} \prod_{k=0}^{b-1} 
\mathbb{E}(\epsilon_0^{\beta_k})M_{\beta_k}^{(n)}=\sum_{\substack{ \beta\in
S_{2p}  \\ \beta_k \mbox{ even}}}\gamma_{ \beta} \prod_{k=0}^{b-1}
M_{\beta_k}.
\end{equation*}
Denote by $L$ the right hand side of the above relation and define $
L'=(b^p-b)^{-1} L$. By using (\ref{H}) we deduce from the previous
lines that

\begin{equation}  \label{rel}
M^{(n+1)}_{2p}=r_n^{2p}b^{1-p}M^{(n)}_{2p}+b^{-p}L+o(1).
\end{equation}
Then by using that $r_n\to 1$ as $n\to\infty$ and the relation $L'=
b^{1-p}L'+b^{-p} L$ we obtain 
\begin{equation*}
M^{(n+1)}_{2p}-L'= r_n^{2p}b^{1-p}(M^{(n)}_{2p}-L')+o(1).
\end{equation*}
This yields both $(\mathcal{P}_{2p})$ and (3) since  $r_n^{2p}b^{1-p}\sim b^{1-p}<1$ as $n\to\infty$.

\medskip

Now suppose that $H<1/2$. Almost the same arguments as when $H=1/2$ yield
the conclusion. The only change is that we have to perform one more
computation to obtain a relation equivalent to (\ref{rel}). Due to the
expression of $r^{(H) }_n$, we have
\begin{eqnarray*}
M^{(n+1)}_{2p}&=&b^{-(1-2H)p}\left (b^{(1-2pH)}M_{2p}^{(n)}+b^{-2pH} L\right
)+o(1) \\
&=& b^{1-p}M_{2p}^{(n)}+b^{-p} L+o(1).
\end{eqnarray*}
\end{proof}

\begin{proposition}
\label{prop2} Let $H\in (-\infty,1/2]$. Let $(W(t))_{t\in [0,1]}$ be a standard Brownian motion. For every $p\ge 1$, the
probability distribution $\mathcal{L}\bigg(\big (B^H_n(t)/a_n^H\big )_{t\in 
\mathcal{D}_p}\bigg )$ converges to $\mathcal{L}\big ((W(t))_{t\in\mathcal{D}
_p}\big )$ as $n\to\infty$.
\end{proposition}

\begin{proof}
Let $p\ge 1$ and denote by $0=t_0<t_1\dots< t_{2^p}=1$ the elements of $
\mathcal{D}_p$. Also, simply denote $B^H_n(t)/a_n^H$ by $\mathcal{Y}_n(t)$
and by $\phi_n$ the characteristic function of $\mathcal{Y}_n(1)$ ($\mathcal{Y}
_n(1)$ is nothing but the random variable $Y_n$ studed in Proposition~\ref{prop1}). By using the fact that in (\ref{fonc0}) the fonctions $
B^H_{n-|w|}(w)$ are mutually independent and identically distributed with $
a_{n-|w|}^H \mathcal{Y}_{n-|w|}(1)$, and also independent of the products $
\epsilon(w|1) \cdots \epsilon (w|n) b^{-nH}$, we can get that for $
(u_w)_{w\in \Sigma_p}\in \mathbb{R}^{2^p}$ and $n>p$ 
\begin{equation*}
\mathbb{E}\left ( e^{i\sum_{w\in\Sigma_p}u_w \big (\mathcal{Y}_n
(t_w+b^{-p})-\mathcal{Y}_n(t_w)\big )}\right )= \mathbb{E}\prod_{w\in
\Sigma_p}\phi_{n-p}\left (u_w b^{-pH} \frac{a_{n-p}^H}{a_n^H}
\prod_{k=1}^p\epsilon(w|k) \right ).
\end{equation*}
It follows from Proposition~\ref{prop1} that $\phi_{n-p}(t)$ goes to $
e^{-t^2/2}$ as $n$ goes to $\infty$. Moreover, $b^{-pH} a_{n-p}^H/a_n^H$
tends to $b^{-p/2}$ as $n$ goes to $\infty$. Thus, applying the dominated
convergence theorem yields 
\begin{eqnarray*}
\lim_{n\to\infty} \mathbb{E}\left ( e^{i\sum_{w\in\Sigma_p}u_w \big (
\mathcal{Y}_n (t_w+b^{-p})-\mathcal{Y}_n(t_w)\big )}\right )&=&\mathbb{E}
\left (\prod_{w\in \Sigma_p} \exp\Big (2^{-1}u_w^2
b^{-p}\prod_{k=1}^p\epsilon(w|k) ^2\Big)\right ) \\
&=& \prod_{w\in \Sigma_p}\exp\big( u_w^2 b^{-p}/2\big )
\end{eqnarray*}
since the $\epsilon(w|k)$ take values in $\{-1,1\}$. This yields the
conclusion.
\end{proof}

\begin{proposition}
\label{prop3} Let $H\in (-\infty,1/2]$. The sequence $\big (\mathcal{L}
(B^H_n/a_n^H)\big )_{n\ge 1}$ of probability distributions on $C([0,1])$ is
tight.
\end{proposition}

\begin{proof}
Let us denote by $\mathcal{Y}_n$ the process $B^H_n/a_n^H$ as in the proof of Proposition~\ref{prop2}. By Theorem 7.3
of \cite{Bil}, since $\mathcal{Y}_n(0)=0$ almost surely for all $n\ge 1$, it
is enough to show that for each positive $\varepsilon$

\begin{equation}  \label{tighness}
\lim_{\delta\to 0}\limsup_{n\to\infty} \mathbb{P}\big (\omega(\mathcal{Y}
_n,\delta)\ge \varepsilon\big )=0
\end{equation}
(the modulus of continuity $\omega(f,\cdot)$ of a continuous function $f$ is defined is Section~\ref{def}).

Fix $H^{\prime }\in (0,1/2)$ and $K$ a positive integer such that $
2K(1/2-H^{\prime })>1$. It follows from the proof of Proposition~\ref{prop1}
that the sequence $\left (\mathbb{E}\left( \mathcal{Y}_{m}(1)^{2K}\right
)
\right )_{m\ge 1}$ is bounded by a constant $C_K$. Moreover, by construction
there exists a constant $C_{b,H}$ such that for all $n>p\ge 1$ we have $
\frac{a_{n-p}^H}{a_n^H}\le C_{b,H} b^{p(H-1/2)}$. By using (\ref{fonc0}) and
a Markov inequality we can get that for $n\ge 2$, $1\le p<n$ and $0\le k\le
b^p-1$ 
\begin{eqnarray*}
&& \mathbb{P}\left ( \left |\mathcal{Y}_n\big ((k+1)b^{-p}\big)- \mathcal{Y}_n\big (kb^{-p}\big )
\right |>b^{-pH^{\prime }}\right ) \\
& \le & b^{2KpH^{\prime }}\mathbb{E}\left (\left |\mathcal{Y}_n\big ((k+1)b^{-p}\big)-
\mathcal{Y}_n\big (kb^{-p}\big )\right |^{2K}\right ) \\
&\le & b^{2K(H^{\prime }-H)p} \left (\frac{a_{n-p}^H}{a_n^H}\right )^{2K} 
\mathbb{E}(\mathcal{Y}_{n-p}(1)^{2K}) \\
&\le & C_KC_{b,H}b^{2K(H^{\prime }-1/2)p}.
\end{eqnarray*}
Now let $\alpha_p=C_K C_{b,H}b^{p(1+2K(H^{\prime }-1/2))}$. By our choice of $
H^{\prime }$ and $K$ the series $\sum_{p\ge 1}\alpha_p$ converge. Moreover,
since for $1\le p<n$  the $b$-adic increments of generation $p$ of $\mathcal{
Y}_n$ have the same probability distribution, we have 
\begin{equation*}
\mathbb{P}\left (\exists \ 0\le k< b^{p},\ \left |\mathcal{Y}_n\big (
(k+1)b^{-p}\big)- \mathcal{Y}_n\big (kb^{-p}\big )\right |>b^{-pH^{\prime
}}\right )\le \alpha_p.
\end{equation*}
On the other hand, if $p\ge n$ and $0\le k< b^{p}$, by construction since
there exists a constant $c_{b,H}<1$ such that $a_n^H\ge c_{b,H}b^{n(1/2-H)}$
we have 
\begin{multline*}
\left |\mathcal{Y}_n\big ((k+1)b^{-p}\big)-\mathcal{Y}_n\big (kb^{-p}\big )
\right |=\frac{b^{-n(H-1)}b^{-p}}{a_n^H} \le \frac{b^{n/2-p}}{c_{b,H}} 
\le \frac{b^{-p/2}}{c_{b,H}}
\le \frac{b^{-pH^{\prime }}}{c_{b,H}}.
\end{multline*}
Let $A_p$ denote the rest $\sum_{j\ge p}\alpha_j$. We deduce from the
previous lines that for all $p\ge 1$, 
\begin{equation*}
\sup_{n\ge 2} \mathbb{P}\left (\exists\ j\ge p,\ \exists \ 0\le k< b^{-j},\
\left |\mathcal{Y}_n\big ((k+1)b^{-j}\big)- \mathcal{Y}_n\big (kb^{-j}\big )
\right |>c_{b,H}^{-1} b^{-jH^{\prime }}\right )\le A_{p}.
\end{equation*}
The event $\left \{\forall\ j\ge p,\ \forall \ 0\le k< b^{-j},\ \left |
\mathcal{Y}_n\big ((k+1)b^{-j}\big)- \mathcal{Y}_n\big (kb^{-j}\big )
\right
|\le c_{b,H}^{-1}b^{-jH^{\prime }}\right\}$ is denoted by $E^n_p$ .
One has $\mathbb{P}(E_p^n)\ge 1-A_p$. A simple adaptation of the proof of
the Kolmogorov-Centsov theorem \cite{Centsov} (see the proof of Proposition~
\ref{prop4}(3) in the next section) shows that on $E^n_p$, we have 
\begin{equation*}
\displaystyle \sup_{n\ge 2}\sup_{\substack{ 0\le s<t\le 1  \\ t-s\le b^{-p}}}
\frac{\left |\mathcal{Y}_n(t)- \mathcal{Y}_n(s)\right |}{(t-s)^{H^{\prime }}}
\le \frac{2(b-1)c_{b,H}^{-1}}{1-b^{-H^{\prime }}}.
\end{equation*}
Consequently, for all $n\ge 2$ we have $\omega\big (\mathcal{Y}_n,b^{-p})\big )\le \frac{
2(b-1)c_{b,H}^{-1}b^{-pH^{\prime }}}{1-b^{-H^{\prime }}}$. This yields 
\begin{equation*}
\inf_{n\ge 2} \mathbb{P}\left (\omega\big (\mathcal{Y}_n,b^{-p}\big )\le 
\frac{2(b-1)c_{b,H}^{-1}b^{-pH^{\prime }}}{1-b^{-H^{\prime }}}\right )\ge
\inf_{n\ge 1} \mathbb{P}(E^n_p)\ge 1-A_p.
\end{equation*}
Since $\lim_{p\to\infty}A_p=0$, the previous inequality yields (\ref
{tighness}).
\end{proof}

\noindent \textit{Proof of Theorem~\ref{th1}.} We use the notations of the
three previous propositions. Suppose that $(\mathcal{Y}_{n_k})_{k\ge 1}$ is
subsequence of $(\mathcal{Y}_n)_{n\ge 1}$ which converges weakly to a
probability distribution $\mathcal{W}_{\infty}$. Due to Proposition~\ref
{prop2}, a process $\mathcal{Y}$ such that $\mathcal{L}(\mathcal{Y})=
\mathcal{W}_{\infty}$ has continuous path and is such that for all $p\ge 1$, 
$\mathcal{L}\big ((\mathcal{Y}(t))_{t\in\mathcal{D}_p}\big )=\mathcal{L}
\big
((\mathcal{W}(t))_{t\in\mathcal{D}_p}\big )$. Since $\bigcup_{p\ge 1}
\mathcal{D}_p$ is dense in $[0,1]$ and we know that the almost sure limit of a
sequence of centered Gaussian variables is a centered Gaussian variable with
variance equal to the limit of the variances, we conclude
that $\mathcal{W}_{\infty}=\mathcal{W}.$ Now the final conclusion comes from
Proposition~\ref{prop3}.

\section{Proof of Theorem~\ref{th2} concerning strong convergence when $1/2<H<1
$}\label{strong}

We first construct in Proposition~\ref{prop5} a stochastic process thanks to
the almost sure pointwise convergence of $B^H_n$ over the set $b$-adic
numbers. We establish regularity properties for this process and then
identify this process as the almost sure uniform limit of $B^H_n$
(Proposition~\ref{identification}) by using a result on vector martingales.
At last we prove the result concerning the Hausdorff and box dimensions of the graph
of the limit $B^H$ of $B^H_n$.

\begin{proposition}
\label{prop5} Let $H\in (1/2,1]$. With probability one

\begin{enumerate}
\item for every $b$-adic number $t$ in $[0,1]$ the sequence $B^H_n(t)$
converges to a limit denoted $B^H_\infty(t)$.

\item The function $B^H_\infty$ defined on the $b$-adic numbers possesses a
(necessarily unique) continuous extension to $[0,1]$ also denoted $B^H_\infty$.

\item The function $B^H_\infty$ belongs to $C^{H^{\prime }}([0,1])$ for all $
H^{\prime }<H$.

\item The pointwise H\"older exponent of $B^H_\infty$ at every point of $
[0,1]$ is equal to $H$.
\end{enumerate}
\end{proposition}

We first establish the following useful result on the martingale $
(B^H_n(1))_{n\ge 1}$.

\begin{lemma}\label{lemmoments}
\label{prop4} Let $H\in (1/2,1]$. The martingale $\big (B^H_n(1)\big )_{n\ge
1}$ is bounded in $L^q$ norm for all $q\ge 1$.
\end{lemma}

\begin{proof}
Denote $B^H_n(1)$ by $Z_n$ as in the proof of Proposition~\ref{prop1}. Since 
$(Z_n)_{n\ge 1}$ is a martingale, the sequence $(\mathbb{E}(Z_n^2))_{n\ge 1}$
is non-decreasing. Consequently, it follows from (\ref{fonc2}), (\ref{m2})
and the fact that $b^{1-2H}<1$ since $H>1/2$ that $(Z_n)_{n\ge 1}$ is
bounded in $L^2$ norm and thus it converges almost surely to a limit $
Z_\infty$. Then, the relation (\ref{Mq}) as well as arguments very similar
to those used in the proof of Proposition~\ref{prop1} show that the sequence 
$\mathbb{E}(Z_n^q)$ converges for every integer $q\ge 1$ as $n$ goes to $
\infty$. In particular it is bounded in $L^{2q}$ for every integer $q\ge 1$.
This implies that $\mathbb{E}(Z_\infty^{2q})<\infty$ for every integer $q\ge
1 $ by the Fatou lemma.
\end{proof}

\begin{proof}
\textit{(of Proposition~\ref{prop5})} (1) Since $B^H_n(0)=0$ for all $n\ge 1$
almost surely, it is enough to establish that for every $p\ge 1$ and $0\le
k< b^{-p}$ the sequence $\big (\Delta^H_n(p,k)\big)_{n\ge 1}$ defined by $
\Delta^H_n(p,k)=B^H_n\big ((k+1)b^{-p}\big )-B^H_n\big (kb^{-p}\big )$
converges almost surely. Indeed, since the set of $b$-adic numbers is
countable, this will imply that with probability one, $\big (\Delta^H_n(p,l)
\big)_{n\ge 1}$ converges for every, $p\ge 1$ and $0\le l<b^p$ as $n$ goes
to $\infty$, thus $B^H_n\big (kb^{-p})=\sum_{l=0}^{k-1}\Delta^H_n(p,l)$
converges for $p\ge 1$ and $1\le k\le b^p$ as $n$ goes to $\infty$.

Now, it is sufficient to notice that given $p\ge 1$, $0\le k< b^{-p}$, and $
w=w_1\cdots w_p$ so that $kb^{-p}=\sum_{j=1}^p w_jb^{-j}$, the relation (\ref
{fonc0}) yields for $n\ge p+1$ 
\begin{equation}  \label{delta}
\Delta^H_n(p,k)= \epsilon(w_1) \cdots \epsilon (w_1\cdots w_p) b^{-pH}
B^H_{n-p}(w)(1).
\end{equation}
The convergence of $\Delta^H_n(p,k)$ then comes from Lemma~\ref{prop4} which
ensures that the martingale $(B^H_{n-p}(w)(1))_{n\ge 1}$ converges to a
limit $B^H_\infty(w)(1)$ since it is bounded in $L^2$-norm.

Let $\Delta^H_\infty(p,k)$ and $B^H_\infty(kb^{-p})$ denote the limit of $
\Delta^H_n(p,k)$ and $B^H_n(kb^{-p})$ respectively. By construction, given $
p\ge 1$, $0\le k< b^{-p}$, and $w=w_1\cdots w_p$ so that $
kb^{-p}=\sum_{j=1}^p w_jb^{-j}$, we have 

\begin{eqnarray} 
\label{delta'}\Delta^H_\infty(p,k)&=& B^H_\infty\big ((k+1)b^{-p}\big )-B^H_\infty\big (
kb^{-p}\big ) \\
 \label{delta"}&= &\epsilon(w_1) \cdots \epsilon (w_1\cdots w_p) b^{-pH} B^H_\infty(w)(1).
\end{eqnarray}

\medskip

\noindent (2) and (3) We adapt the proof of the Kolmogorov-Centsov theorem 
\cite{Centsov,KS} which uses the dyadic basis while we work in any basis $b$. Let $H^{\prime }\in (0,H)$. Fix an integer $
K>1/2(H-H^{\prime })$. Due to (\ref{delta'}) and (\ref{delta"}), for $p\ge 1$
we have

\begin{equation*}
\alpha_p:=\mathbb{P}\left (\exists \ 0\le k<b^p,\ |\Delta^H_\infty(p,k)|\ge
b^{-pH^{\prime }}\right )\le b^{p(1+2K(H^{\prime }-H)} \mathbb{E}\big (
B^H_\infty(1)^{2K}\big).
\end{equation*}
Since $\sum_{p\ge 1}\alpha_p<\infty$, the Borel-Cantelli implies that with
probability 1, there exists $n_0$ such that 
\begin{equation}  \label{Centsov0}
\sup_{0\le k<b^n}|\Delta^H_\infty(n,k)|< b^{-nH^{\prime }}, \quad \forall\
n\ge n_0 .
\end{equation}
Now we fix $n\ge n_0$ and show that for all $m>n$,

\begin{equation}  \label{Centsov}
|B^H_\infty(t)-B^H_\infty(s)|\le 2(b-1)\sum_{j=n+1}^m b^{-H^{\prime
}j},\quad \forall\ t,s\in \mathcal{D}_m,\ 0<t-s<b^{-n}.
\end{equation}
If $m=n+1$, one has $s=kb^{-(n+1)}$ and $t=k^{\prime -(n+1)}$ with $
0<k^{\prime }-k<b$, so due to (\ref{Centsov0}) we have $|B^H_\infty(t)-B^H_
\infty(s)|\le (k' -k)b^{-(n+1)H'}$, hence the conclusion.

Suppose that (\ref{Centsov}) holds for $n+1\le m \le M-1$. Let $t,s\in 
\mathcal{D}_M$ such that $0<t-s<b^{-n}$ and consider $t_1=\max\{u\in 
\mathcal{D}_{M-1}: u\le t\}$ and $s_1=\min\{u\in \mathcal{D}_{M-1}: u\ge s\}$
. One has $s\le s_1 \le t_1\le t$, $t_1-s_1<b^{-n}$, $s_1-s\le (b-1) b^{-M}$
and $t-t_1\le (b-1)b^{-M}$. Now, since $s_1$ and $t_1$ belong to $\mathcal{D}
_{M-1}\subset \mathcal{D}_M$, property (\ref{Centsov0}) implies that $
|B^H_\infty(s)-B^H_\infty(s_1)|\le (b-1) b^{-MH^{\prime }}$ and $
|B^H_\infty(t)-B^H_\infty(t_1)|\le (b-1) b^{-MH^{\prime }}$. Moreover, since
(\ref{Centsov}) holds for $m=M-1$ one has $|B^H_\infty(t_1)-B^H_\infty(s_1)|
\le 2(b-1)\sum_{j=n+1}^{M-1}b^{-H^{\prime }j}$. This is enough to get (\ref
{Centsov}) for $m=M$.

\medskip

Property (\ref{Centsov}) being established for all $n\ge n_0$, taking $t,s\in
\mathcal{D}$ such that $0<|t-s|<b^{-n_0}$ and $n$ the integer such that $
b^{-(n+1)}\le |t-s|<b^{-n}$, since both $t$ and $s$ belong to $\bigcup_{p>n}
\mathcal{D}_p$ we deduce from (\ref{Centsov}) that 
\begin{equation*}
|B^H_\infty(t)-B^H_\infty(s)|\le 2(b-1)\sum_{j=n+1}^\infty b^{-H^{\prime
}j}\le \frac{2(b-1)}{1-b^{-H^{\prime }}}|t-s|^{H^{\prime }}.
\end{equation*}
This is enough to construct on $[0,1]$ a unique continuous extension of $
B^H_\infty$. As a consequence of what preceeds, this extension belongs to $
C^{H^{\prime }}([0,1])$ for all $H^{\prime }<H$.

\medskip

\noindent (4) We need the following lemma which describes the asymptotic
behavior of the characteristic function of $B^H_\infty(1)$. This lemma will
be also useful in finding a lower bound for the Hausdorff dimension of the
graph of $B^H$.

\begin{lemma}
\label{fourier} Let $\varphi$ stand for the characteristic function of $
B^H_\infty(1)$. There exists $\rho\in(0,1)$ such that $\varphi(t)=O\big (
\rho^{|t|^{1/H}}\big )\ \ (|t|\to\infty)$. In particular, the probability
distribution of $B^H_\infty(1)$ possesses an infinitely differentiable
density.
\end{lemma}

\begin{proof}
Since $\mathbb{E}(B^H_\infty(1))=1$, the probability distribution of $
B^H_\infty(1)$ is not concentrated at 0 and thus for every $\eta>0$ there
exists $\alpha\in (0,\eta)$ and $\gamma<1$ such that $\sup_{t, |t|\in
[\alpha,b^H\alpha]}|\varphi(t)|\le \gamma$.

Now, using the fact that 
\begin{equation*}
\varphi(t)=\left [p_{b,H}^+\varphi\big (b^{-H}t\big )+p_{b,H}^-\varphi\big (-b^{-H}t\big )\right
]^b,
\end{equation*}
one obtains by induction that 
\begin{equation*}
\sup_{t,\ |t|\in [b^{kH}\alpha,b^{(k+1)H]
}\alpha]}|\varphi(t)|\le \gamma^{b^k}\quad (\forall\ k\ge 0).
\end{equation*}
Since $|t|^{1/H}\le b\alpha^{1/H} b^{k}$ for $|t|\in [b^{kH}\alpha,b^{(k+1)H}\alpha]$, the
conclusion follows with $\rho=\gamma^{1/ b\alpha^{1/H}}$.

The rate of decay of $\varphi$ at $\infty$ yields the conclusion regarding
the probability distribution of $B^H_\infty(1)$.
\end{proof}

It follows from Lemma~\ref{fourier} that $\mathbb{E}(|B^H_\infty(1)|^{-
\gamma})<\infty$ for all $\gamma\in (0,1)$. This will be used with $
\gamma=1/2$ in what follows.

We next use an approach similar to that used for the study of the pointwise
H\"older exponents of  Brownian motion \cite{Erdos,KS}.

Let $\varepsilon>0$. We show that the subset $\mathcal{O}$ of $\Omega$ of
points $\omega$ such that the corresponding path $B^H_\infty$ possesses
points at which the pointwise H\"older exponent is at least $H+\varepsilon$
is included in a set of null probability.

We fix an integer $K>4/\varepsilon $ and denote by $n_K$ the smallest
integer $n$ such that $Kb^{-n}\le 1$. For $t\in [0,1]$ and $n\ge n_K$ , consider
$S^K_n(t)$ a subset of $[0,1]$ made of $K+1$ consecutive $b$-adic
numbers of generation $n$ such that $t\in [\min\, S^K_n(t),\max\, S^K_n(t)]$. 
Also denote by $\boldsymbol{S}^K_n(t)$ the set of $K$ consecutive $b$-adic
intervals delimited by the elements of $S^K_n(t)$. If the pointwise H\"older
exponent at $t$ is larger than or equal to $H+\varepsilon$ then for $n$
large enough one has necessarily $\sup_{s\in
S^K_n(t)}|B_\infty^H(s)-B_\infty^H(t)|\le (Kb^{-n})^{H+\varepsilon/2}$, so
that $\sup_{I\in \boldsymbol{S}^K_n(t)} |\Delta B_\infty^H(I)|\le 2
(Kb^{-n})^{H+\varepsilon/2}$.

Know let $\boldsymbol{S}^K_n$ be the set made of all $K$-uple of consecutive $
b$-adic intervals of generation $n$, and if $S\in \boldsymbol{S}^K_n$, denote the event $\big \{ \sup_{I\in S} |\Delta
B_\infty^H(I)|\le 2 (Kb^{-n})^{H+\varepsilon/2}\big \}$ by $E_{S}$. The previous lines
show that 
\begin{equation*}
\mathcal{O}\subset \mathcal{O}^{\prime }=\bigcap_{n\ge n_K}\bigcup_{p\ge
n}\bigcup_{S\in \boldsymbol{S}^K_p}E_{S}.
\end{equation*}
By construction, if $S\in \boldsymbol{S}^K_p$, $\big (|\Delta
B_\infty^H(I)|\big)_{I\in S}$ is equal to $(
b^{-pH}Y_I)_{I\in S}$, where the $K$ random variables $Y_I$ are
mutually independent and identically distributed with $|B^H_\infty(1)|$.
Consequently, $\mathbb{P}(E_{S})$ depends only on $K$ and $p$ and 
\begin{equation*}
\mathbb{P}(E_{S})\le \left [ \mathbb{P}(|B^H_\infty|\le 2
(Kb^{-p})^{\varepsilon/2})\right ]^K\le \sqrt{2}^K K^{K \varepsilon/4} b^{-p
K \varepsilon/4}\left [\mathbb{E}(|B^H_\infty(1)|^{-1/2})\right ]^K.
\end{equation*}
Since the cardinality of $\boldsymbol{S}^K_p$ is less than $b^p$, this yields $
\mathbb{P}\big(\bigcup_{S\in \boldsymbol{S}^K_p}E_{S}\big )\le C b^p
b^{-p K \varepsilon/4}$, with $C=\sqrt{2}^K K^{K \varepsilon/4}\left [
\mathbb{E}(|B^H_\infty(1)|^{-1/2})\right ]^K$. Due to our choice for $K$,
this implies that the series $\sum \mathbb{P}\big(\bigcup_{S\in \boldsymbol{S}
^K_p}E_{S}\big )$ converges and $\mathbb{P}(\mathcal{O}^{\prime
})=0$.
\end{proof}

The next proposition makes it possible to conclude that the random sequence
of functions $(B^H_n)_{n\ge 1}$ converges almost surely uniformly to the
function $B^H_\infty$ constructed previously. The same kind of approach is
used to establish the convergence of continuous function-valued martingales
related to multiplicative processes on a homogeneous or Galton Watson tree
in \cite{Joffe,Biggins,Barral}, but the context in the mentioned papers is
rather different from the present one because the martingales $M_n(s)$
considered there take the form $\sum_{w\in \Sigma_n} \prod_{k=1}^nW(w|k)(s)$
where the random weights $W(w)(s)$ depend smoothly on the parameter $s$
belonging to some open subset of $\mathbb{R}^d$ independently of $\Sigma$ (and more generally a super-critical Galton-Watson tree),
while for $B^H_n(s)$ the parameter $s$ is a generic point in $\Sigma$
(identified with $[0,1]$).

\begin{proposition}
\label{identification} One has $\mathbb{E}\left
(\|B^H_\infty\|_\infty
\right
)<\infty$. Consequently, with probability 1, $B^H_n$ converges
uniformly to $B^H_\infty$.
\end{proposition}

\begin{proof}
In fact we are going to prove that $\mathbb{E}\left
(\|B^H_\infty\|_
\infty^2\right )<\infty$. Define 
\begin{equation*}
\widetilde Z_n=\sup_{t\in\bigcup_{p\ge 1}\mathcal{D}_p}|B^H_n(t)|, \ \widetilde Z_{n}(k)=
\sup_{t\in\bigcup_{p\ge 1}\mathcal{D}_p}|B^H_n(k)(t)|,\ 0\le k\le b-1.
\end{equation*}
Due to (\ref{fonc}) one has 
\begin{equation*}
\widetilde Z_n\le \max_{0\le k\le b-1}b^{-H}\widetilde Z_{n-1}(k)+|B^H_n(kb^{-1})|\le \sum_{k=0}^{b-1}b^{-H}\widetilde Z_{n-1}(k)+|B^H_n(kb^{-1})|
\end{equation*}
Thus 
\begin{eqnarray*}
\mathbb{E}(\widetilde Z_n^2)&\le& \sum_{k=0}^{b-1}\mathbb{E}\left
(b^{-2H}\widetilde Z_{n-1}(k)^2+2\widetilde Z_{n-1}(k)|B^H_n(kb^{-1})|+|B^H_n(kb^{-1})|^2\right )
\\
&=&b^{1-2H}\mathbb{E}(\widetilde Z_{n-1}^2)+ 2\mathbb{E}\left
(\sum_{k=0}^{b-1}\widetilde Z_{n-1} (k)|B^H_n(kb^{-1})|\right )+\sum_{k=0}^{b-1}\mathbb{E}
(|B^H_n(kb^{-1})|^2 ) \\
&\le & b^{1-2H}\mathbb{E}(\widetilde Z_{n-1}^2)+2\mathbb{E}(\widetilde Z_{n-1}^2)^{1/2} 
\sum_{k=0}^{b-1}\left\|B^H_n(kb^{-1})\right \|_2+
\sum_{k=0}^{b-1}\| B^H_n(kb^{-1})\|_2^2.
\end{eqnarray*}
Now we use the fact that the sequence 
$\sup_{0\le k\le b-1}\left\| B^H_n(kb^{-1})\right \|_2$ is bounded due to the proof of Proposition~\ref{prop5}. Thus there exists $C>0$ such that for all 
$n\ge 1$ 
\begin{equation}  \label{infty}
\mathbb{E}(\widetilde Z_n^2)\le f\big (\mathbb{E}(\widetilde Z_{n-1}^2)\big ),\ \mbox{with } f(x)=
b^{1-2H}x+C\sqrt {x}+C.
\end{equation}
Since $b^{1-2H}<1$, there exists $x_0> 0$ such that $f(x)<x$ for all $x>x_0$
. This remark together with (\ref{infty}) implies that $\mathbb{E}(\widetilde Z_n^2)\le
\max \left (x_0, f\big (\mathbb{E}(\widetilde Z_{1}^2)\big )\right )$ for all $n\ge 2$.

To conclude, we use Proposition V-2-6 in \cite{Neveu} which ensures that since $
\mathcal{C}([0,1])$ is a complete separable Banach space and $
\|B^H_\infty\|_{L^1}<\infty$, $B^H_\infty$ is the almost sure limit of $
\widetilde B^H_n=\mathbb{E}\big (B^H_\infty|\sigma(\epsilon (w): |w|\le n)
\big )$. Furthermore, given $n\ge 1$ and $w\in\Sigma_n$, one can show by
induction on $p\ge 0$ that, with probability 1, $\widetilde B^H_n(t_{w\cdot
u})=B^H_n(t_{w\cdot u})$ for all $u\in \Sigma_p$. This implies that $
\widetilde B^H_n=B^H_n$ almost surely since these functions coincide over $\mathcal{D}$.
\end{proof}

\begin{proposition}
Let $H\in (1/2,1]$. With probability 1, the Hausdorff and box dimensions of
the graph of $B_H$ are equal to $2-H$.
\end{proposition}

We shall need an additional notation. If $w\in\Sigma^*$ and $J=\pi ([w])$
then we define $\boldsymbol{\epsilon}(J)=\prod_{k=1}^{|w|}\epsilon(w|k)$.

\begin{proof}
Let us denote by $\Gamma_H$ the graph $\left \{\big (t,B_H(t)\big ): t\in
[0,1]\right \}$ of $B_H$.

At first, the fact that $2-H$ is an upper bound for the box dimension of the
graph of $B_H$ comes from the fact that $B_H\in C^{H^{\prime }}([0,1])$ for
all $H^{\prime }<H$ (see \cite{Falc} Ch. 11).

To find the sharp lower bound $2-H$ for the Hausdorff dimension of $\Gamma_H$
, we use the method consisting in showing that with probability 1, the
measure on this graph obtained as the image of the Lebesgue measure
restricted to $[0,1]$ by the mapping $t\mapsto \big (t,B_H(t)\big )$ has a
finite energy with respect to the Riesz Kernel $u\in\mathbb{R}
^2\setminus\{0\}\mapsto \|u\|^{-\gamma}$ for all $\gamma <2-H$ (see \cite{Falc} Ch. 11
for more details). This property holds if we show that for all $\gamma <2-H$

\begin{equation*}
\int_{[0,1]^2}\mathbb{E}\left (\frac{1}{\sqrt{|t-s|^2+|B_H(t)-B_H(s)|^2}
^\gamma}\right )\ dtds <\infty.
\end{equation*}
If $I$ is a closed subinterval of $[0,1]$, we denote by $\mathcal{G}(I)$ the
set of closed $b$-adic intervals of maximal length included in $I$, and then 
$m_I=\min\bigcup_{J\in\mathcal{G}(I)}J$ and $M_I =\max\bigcup_{J\in\mathcal{G
}(I)}J$.

Let $0<s<t<1$ be two non $b$-adic numbers. We define two sequences $
(s_p)_{p\ge 0}$ and $(t_p)_{p\ge 0}$ as follows. Let $s_0=m_{[s,t]}$ and $
t_0=M_{[s,t]}$. Then let define inductively $(s_p)_{p\ge 1}$ and $
(t_p)_{p\ge 1}$ as follows: $s_p=m_{[s,s_{p-1}]}$ and $t_p=M_{[t_{p-1},t]}$.
Let us denote by $\mathcal{C}$ the collection of intervals containing $[s_0,t_0]$ and the intervals $[s_p,s-{p-1}]$ and $[t_{p-1},t_p]$, $p\ge 1$. Every interval $I\in \mathcal{C}$ is the union of at most $b-1$
intervals of the same generation $n_I$, the elements of $\mathcal{G}(I)$,
and we have $\Delta B_H(I)=\sum_{J\in\mathcal{G}(I)} \boldsymbol{\epsilon}
(J)b^{-n_IH}Y_J$. 

By construction, we have $\min_{I\in\mathcal{C}
}n_I=n_{[s_0,t_0]}$ and $(t-s)/b\le b^{-n_{[s_0,t_0]}}\le (t-s)$.

Also, all the random variables $Y_I$ are mutually independent and
independent of $\mathcal{T}_{\mathcal{C}}=\sigma (\boldsymbol{\epsilon}
(J):J\in \mathcal{G}(I),\ I\in\mathcal{C})$.

Now, we write 
\begin{equation*}
B_H(t)-B_H(s)=b^{-n_{[s_0,t_0]}H}\left (\sum_{J\in \mathcal{G}([s_0,t_0])} 
\boldsymbol{\epsilon} (J)Y_J+ Z(s,s_0)+Z(t_0,t)\right),
\end{equation*}
where 
\begin{equation*}
\begin{cases}
\displaystyle Z(s,s_0)=\lim_{p\to\infty} \sum_{0\le k\le
p}b^{(n_{[s_0,t_0]}-n_{[s_{k+1},s_k]})H} \sum_{J\in\mathcal{G}
([s_{k+1},s_k])} \boldsymbol{\epsilon} (J)Y_J \\ 
\displaystyle Z(t_0,t)=\lim_{p\to\infty} \sum_{0\le k\le
p}b^{(n_{[s_0,t_0]}-n_{[t_k, t_{k+1}]})H} \sum_{J\in\mathcal{G}([t_k,
t_{k+1}])} \boldsymbol{\epsilon} (J)Y_J.
\end{cases}
\end{equation*}
Let $\mathcal{Z}(t,s)=\sum_{J\in \mathcal{G}([s_0,t_0])} \boldsymbol{\epsilon
} (J)Y_J+ Z(s,s_0)+Z(t_0,t)$ and fix $J_0\in \mathcal{G}([s_0,t_0])$.
Conditionally on $\mathcal{T}_{\mathcal{C}}$, $\mathcal{Z}(t,s)$ is the sum
of $\pm Y(J_0)$ plus a random variable $U$ independent of $Y(J_0)$.
Consequently, the probability distribution of $\mathcal{Z}(t,s)$
conditionally on $\mathcal{T}_{\mathcal{C}}$ possesses a density $f_{t,s}$
and $\|\widehat{f_{t,s}}\|_{L^1}\le \|\varphi\|_{L^1}$, where $\varphi$ is
the characteristic function of $Y(J_0)$ studied in Lemma~\ref{fourier}.

Thus, for $\gamma<2-H$ we have

\begin{eqnarray*}
\mathbb{E}\left (\frac{1}{\sqrt{|t-s|^2+|B_H(t)-B_H(s)|^2}^\gamma}|\mathcal{T
}_{\mathcal{C}}\right )&=& \int_{\mathbb{R}}\frac{f_{t,s}(u)}{\sqrt{
|t-s|^2+b^{-2n_{[s_0,t_0]}H}u^2}^\gamma}\, du\\
&\le& \int_{\mathbb{R}}\frac{f_{t,s}(u)}{\sqrt{
|t-s|^2+b^{-2H}(t-s)^{2H}u^2}^\gamma}\, du \\
&=&|t-s|^{1-H-\gamma} \int_{\mathbb{R}}\frac{f_{t,s}(|t-s|^{1-H}v)}{\sqrt{
1+b^{-2H}v^2}^\gamma}\, dv.
\end{eqnarray*}

The function $f_{t,s}$ is bounded independently of $t,\ s$ and $\mathcal{T}_{\mathcal{%
C}}$ since it is bounded by $\|\widehat{f_{t,s}}\|_{L^1}$ and we just saw that
this number is bounded by $\|\varphi\|_{L^1}$. It follows that 
\begin{equation*}
\mathbb{E}\left (\frac{1}{\sqrt{|t-s|^2+|B_H(t)-B_H(s)|^2}^\gamma}\right
)\le \|\varphi\|_{L^1}|t-s|^{1-H-\gamma}\int_{\mathbb{R}}\frac{dv}{\sqrt{1+b^{-2H}v^2}^\gamma}
.
\end{equation*}
This yields the conclusion.
\end{proof}

\section{Proof of Theorems~\ref{th3} and \ref{th4} concerning functional CLT when $1/2<H<1$}
\noindent
{\it Proofs of Theorem~\ref{th3} and ~\ref{th4}.} We proceed in three steps as for Theorem~\ref{th1}.

Let $a_n^H=\sigma_Hb^{n(1/2-H)}$. For $w\in\Sigma^*$ and $n\ge 1$, let $Y_n(w)=\big (B^H(w)-B^H_{n}(w)\big )(1)/a_n^H$, and simply denote $Y_n(\emptyset)$ by $Y_n$. By construction, $\mathcal{L}(Y_n(w))=\mathcal{L}(Y_n)$. Also, for $n\ge 1$ let $\mathcal{X}^H_{n}=B^H-B^H_n$. 

\medskip

Step 1: We leave the reader verify that 
\begin{equation}\label{zzz}
Y_n=\frac{1}{\sqrt{b}}\sum_{k=0}^{b-1} \epsilon (k)Y_{n-1}(k),
\end{equation}
where $Y_{n-1}(k)\sim Y_{n-1}$, and the $Y_{n-1}(k)$'s are centered, mutually independent, and independent of the $\epsilon(k)$'s. 

It is then straightforward to show by induction that properties (1), (2) and (3) of the proof of Proposition~\ref{prop1} hold (with the new sequence $(Y_n)$ considered in the present proof). 

Thus $X^H_n(1)=\mathcal{X}^H_n(1)$ converges weakly to $\mathcal{N}(0,1)$.

\medskip

Step 2: Let $(W(t))_{t\in [0,1]}$ be a standard Brownian motion. For every $p\ge 1$, the
probability distributions $\mathcal{L}\bigg(\big (X^H_{p+n}(t)\big )_{t\in 
\mathcal{D}_p}\bigg )= \mathcal{L}\bigg(\big (\mathcal{X}^H_{p+n}(t)\big )_{t\in 
\mathcal{D}_p}\bigg )$ converge to $\mathcal{L}\big ((W(t))_{t\in\mathcal{D}
_p}\big )$ as $n\to\infty$. This is obtained by following the same approach as in the proof of Proposition~\ref{prop2}, as well as the step 1 and the fact that if $w\in\Sigma_p$, for all $n\ge p$ we have 
\begin{equation}\label{yyy}
\Delta \mathcal{X}^H_n(I_w)=\Delta X^H_n(I_w)=b^{-p/2}Y_{n-p}(w).
\end{equation}
\medskip

Step 3: To see that the sequences $\big (\mathcal{L}(X^H_n)\big )_{n\ge 1}$ and $\big (\mathcal{L}(\mathcal{X}^H_n)\big )_{n\ge 1}$ are tight, we follow the same approach as in Proposition~\ref{prop3}.  

For $n\ge 2$ and $p\ge 1$, we notice that:

\medskip

If $1\le p\le n$ and $w\in \Sigma_p$ then we have $\Delta \mathcal{X}^H_n(I_w)=\Delta X^H_n(I_w)=b^{-p/2}Y_{n-p}(w)$ (this is (\ref{yyy})).

\medskip

If $p>n$ and $w\in \Sigma_p$ then we have 
\begin{eqnarray*}
|\Delta X^H_n(I_w)|= |\Delta X^H_n(I_{w|n})|b^{n-p}&=&\sigma_H^{-1}b^{n/2-p} \big |B^H(w|n)(1)-1\big|\\
&\le &\sigma_H^{-1}b^{-p/2} \big |B^H(w|n)(1)-1\big|.
\end{eqnarray*}
and 
\begin{eqnarray*}
|\Delta \mathcal{X}^H_n(I_w)|&=&(a^H_n)^{-1}\left |\Delta B^H(I_w)-\Delta B^H_n(I_w)\right |\\
&=&\sigma_H^{-1}b^{n(H-1/2)}\left |b^{-pH}B^H(w)(1)\prod_{k=1}^p\epsilon (w|k)-b^{-nH}b^{n-p}\prod_{k=1}^n\epsilon (w|k)\right |
\\
&=&\sigma_H^{-1}b^{-p/2} \left |b^{(p-n)(1/2-H)}B^H(w)(1)\prod_{k=n+1}^p\epsilon (w|k)-b^{-(p-n)/2}\right |.
\end{eqnarray*}
Now, since by Lemma~\ref{lemmoments} and the step 1 the sequences $Y_{n}(w)$, $\big |B^H(w|n)(1)-1\big|$ and $\left |b^{(p-n)(1/2-H)}B^H(w)(1)\prod_{k=n+1}^p\epsilon (w|k)-b^{-(p-n)/2}\right |$ are bounded in $L^q$ for all $q\ge 1$ independently of $w$, we can deduce (by using an approach similar to that used in the proof of  Proposition~\ref{prop3}) from the previous estimates of 
$|\Delta \mathcal{X}^H_n(I_w)|$ and $|\Delta X^H_n(I_w)|$  that for every $H'<1/2$, there exists a positive sequence $(\beta_p)_{p\ge 1}$ such that $\sum_{p\ge 1}\beta_p<\infty$ and 
$$
\sup_{n\ge 2} \mathbb{P}(\exists\ w\in \Sigma_p: |\Delta X^H_n(I_w)|>b^{-pH'})+\sup_{n\ge 2} \mathbb{P}(\exists\ w\in \Sigma_p: |\Delta \mathcal{X}^H_n(I_w)|>b^{-pH'})\le \beta_p.
$$
In view of the proof of Proposition~\ref{prop3}, this is enough to establish the desired tightness.

\end{document}